\documentclass[11pt]{article}
\usepackage[a4paper,left=2cm,right=2cm,top=2.5cm,bottom=2.5cm]{geometry}

\usepackage{mathtools,amsmath,amsfonts}
\usepackage[colorlinks=true,pagebackref=true]{hyperref}

\hypersetup{
    linkcolor=blue,
    citecolor=green,
    urlcolor=blue
}
\usepackage{amsthm, amssymb}
\usepackage{booktabs}
\usepackage{enumitem}

\newtheorem{theorem}{Theorem}[section]
\newtheorem{definition}[theorem]{Definition}
\newtheorem{proposition}[theorem]{Proposition}
\newtheorem{corollary}[theorem]{Corollary}
\newtheorem{lemma}[theorem]{Lemma}
\newtheorem{remark}[theorem]{Remark}
\newtheorem{example}[theorem]{Example}
\newtheorem{examples}[theorem]{Examples}

\numberwithin{equation}{section}

\newcommand{\mult}{\hbox{\rm mult}}

\def\C{{\mathcal C}}

\title{Concentration of plane branches: building a complete equisingularity invariant}
\author{E. R. Garc\'{\i}a Barroso and W.G. Hern\'andez-Yanes}
\date{}

\begin{document}

\maketitle

\begin{abstract}
We study the concentration of the semigroup of values of a plane algebroid branch over an algebraically closed field of arbitrary characteristic. We give an explicit formula for it in terms of the multiplicities of the branch and of its strict transform; in particular, it is an invariant of the equisingularity class. Applying it to the successive strict transforms, we associate a concentration sequence with each branch. This sequence is not a complete equisingularity invariant: we describe all the multiplicity sequences sharing a given concentration sequence, give an algorithm to compute them, and count them in terms of certain positions of the sequence, which we call free. We also characterize the sequences of positive integers that arise as concentration sequences.
Finally, we build a complete equisingularity invariant: the concentration sequence, together with the conductor of the semigroup of values, determines the equisingularity class of the branch, while neither of them does so alone. In characteristic zero, this means that the concentration sequence and the
Milnor number form a complete invariant.


\end{abstract}

\section{Introduction}

Let $K$ be an algebraically closed field of arbitrary characteristic, and let $C$ be a plane algebroid branch, that is, the germ at the origin of the curve defined by an irreducible power series $f\in K[[x,y]]$ with $f(0,0)=0$. Put $C^{(0)}=C$ and, for $i\geq 1$, let $C^{(i)}$ be the \emph{strict transform} of $C^{(i-1)}$ by the blowup of its singular point; we also write $\widetilde C=C^{(1)}$. Iterating blowups we obtain, after finitely many steps, a smooth strict transform, that is, one of multiplicity~$1$; this process is the \emph{resolution of the singularity} of $C$. The multiplicities $\mult(C^{(0)}),\ldots,\mult(C^{(r)})$, where $C^{(r)}$ is the first smooth strict transform, form the \emph{multiplicity sequence} of $C$. Two branches are \emph{equisingular} if they have the same multiplicity sequence or, equivalently, the same semigroup of values (see~\cite{Zariski-book} for $K=\mathbb{C}$ and~\cite{Campillo} for arbitrary characteristic).

\medskip

We recall now some notions on numerical semigroups. A \emph{numerical semigroup} is a subset $S$ of $\mathbb{N}$ that is closed under addition, contains $0$, and has finite complement in $\mathbb{N}$. Every numerical semigroup is finitely generated, and it has a unique \emph{minimal system of generators}, denoted by $\mathrm{msg}(S)=\{v_0<\dots<v_g\}$, which satisfies $\gcd(\mathrm{msg}(S))=1$ (see~\cite{R-G2009}). The smallest positive element of $S$ is its \emph{multiplicity} $m(S)$, and the elements of $\mathbb{N}\setminus S$ are its \emph{gaps}. The largest gap is the \emph{Frobenius number} $F(S)$, and $c(S)=F(S)+1$ is the \emph{conductor} of $S$.

\medskip

Associated with a branch $C$ we have its \emph{semigroup of values}:
\[
\Gamma(C):=\{i_0(f,h) : h\in K[[x,y]],\ f \text{ does not divide } h\},
\]
where $i_0(f,h):=\dim_K K[[x,y]]/(f,h)$ is the \emph{intersection multiplicity} of $C$ with the germ defined by $h(x,y)=0$. It is well known that $\Gamma(C)$ is a numerical semigroup whose multiplicity is $\mult(C)$.

\medskip

The \emph{concentration} of a numerical semigroup $S$ is $\C(S)=\max\{\mathrm{next}_S(s)-s : s\in S\setminus\{0\}\}$, where $\mathrm{next}_S(s)=\min\{x\in S : x>s\}$. This notion was introduced in~\cite{R-B-T2022}, where the numerical semigroups with concentration~$2$ were characterized, and these results were extended in~\cite{R-B-T2025} to arbitrary concentration.

\medskip

In this paper we study the concentration of the semigroup of values $\Gamma(C)$ of a plane algebroid branch $C$ and of the semigroups of values of its strict transforms, and we use it to build a complete equisingularity invariant. In Section~\ref{sec:preliminares} we bound the concentration of a numerical semigroup in terms of its first two generators.

\medskip

In Section~\ref{sec:comun}, Proposition~\ref{proposicion concentracion curva} shows that $\C(\Gamma(C))=\max\{\mult(\widetilde C),\mult(C)-\mult(\widetilde C)\}$. Hence $\C(\Gamma(C))$ depends only on the first two terms of the multiplicity sequence of $C$, and the \emph{concentration sequence}
\[
\mathcal{S}=\bigl(\C(\Gamma(C^{(0)})),\ldots, \C(\Gamma(C^{(r-1)}))\bigr)
\]
is an invariant of the equisingularity class of $C$.

\medskip

In Section~\ref{sec:estructura} we study the set of multiplicity sequences of the branches with concentration sequence $\mathcal{S}$. Theorem~\ref{teo:descripcion} characterizes them by local conditions on consecutive terms, which gives an algorithm that computes all of them. Theorem~\ref{teo:producto} describes the structure of this set: two of its elements differ only in their first three terms and at certain positions, which we call free and which can be read off from $\mathcal{S}$, and the choices at these positions are independent. As a consequence, the number of equisingularity classes of branches with concentration sequence $\mathcal{S}$ is $t\,2^{f}$, where $t\in\{1,2,3,5\}$ and $f$ is the number of free positions. In particular, $\mathcal{S}$ is not a complete invariant, and Corollary~\ref{cor:cota-r} bounds this number by $5\cdot 2^{\lfloor (r-2)/3\rfloor}$. We finish the section with Proposition~\ref{prop:imagen}, which characterizes the sequences of positive integers that are concentration sequences of plane branches, by conditions read from right to left.

\medskip

In Section~\ref{sec:extendida} we show that adding the conductor $c(\Gamma(C))$ to $\mathcal{S}$ gives a complete invariant of the equisingularity class: the \emph{extended concentration sequence} $\mathcal{S}^e:=(c(\Gamma(C));\mathcal{S})$ determines the multiplicity sequence of $C$, and hence $\Gamma(C)$ (Theorem~\ref{teo:invariante-completo}). In characteristic zero, the conductor can be replaced by the Milnor number (see Remark~\ref{rem:milnor}).

\section{Concentration of numerical semigroups}\label{sec:preliminares}

In this section we collect the results on numerical semigroups that we need. Recall that the concentration of a numerical semigroup $S$ is $\C(S)=\max\{\mathrm{next}_S(s)-s : s\in S\setminus\{0\}\}$, where $\mathrm{next}_S(s)=\min\{x\in S : x>s\}$. We start with a bound for the concentration in terms of the multiplicity, which is observed without proof in~\cite[p.~63]{R-B-T2025}, and a sufficient condition for equality; we include its short proof for completeness.

\begin{proposition}\label{proposicion la multiplicidad acota la concentracion}
    Let $S$ be a numerical semigroup and $m$ its multiplicity. Then, $\C(S)\leq m$. Moreover, if $\mathrm{msg}(S)=\{v_0<v_1<\cdots<v_g\}$ with $g\geq 1$ and $v_1>2v_0$, then $\C(S)=m$.
\end{proposition}

\begin{proof}
Let $s\in S\setminus\{0\}$. Then $s+m\in S$ and $s+m>s$, so $\mathrm{next}_S(s)\leq s+m$, that is, $\mathrm{next}_S(s)-s\leq m$. Hence $\C(S)\leq m$. Now assume that $g\geq 1$ and $v_1>2v_0$, where $v_0=m$. An element $x\in S$ with $m<x<2m$ is not the sum of two nonzero elements of $S$, so it is a minimal generator different from $v_0$, and hence $x\geq v_1>2m$, a contradiction. Thus, as $2m\in S$, we get $\mathrm{next}_S(m)=2m$, so $\C(S)\geq m$, and therefore $\C(S)=m$.
\end{proof}

To compute the concentration, it is enough to look at the minimal
generators:

\begin{lemma}[{\cite[Corollary~3]{R-B-T2025}}]\label{corolario concentracion generadores minimales}
    Let $S$ be a numerical semigroup. Then, \[\C(S)=\max\{\mathrm{next}_S(x)-x : x\in\mathrm{msg}(S)\}.\]
\end{lemma}

When the first two generators satisfy $v_1<2v_0$, they give a better bound than the multiplicity.

\begin{proposition}\label{prop:cota concentracion}
Let $S$ be a numerical semigroup with $\mathrm{msg}(S)=\{v_0<v_1<\cdots<v_g\}$, $g\geq 1$, and $v_1<2v_0$. Then \begin{equation}\label{ineq:concentration} \C(S)\leq\max\{v_1-v_0,\;2v_0-v_1\}.
\end{equation}
\end{proposition}

\begin{proof}
Let $e_1=\gcd(v_0,v_1)$ and $S'=\left\langle\frac{v_0}{e_1},\frac{v_1}{e_1}\right\rangle$. Since $\gcd\left(\frac{v_0}{e_1},\frac{v_1}{e_1}\right)=1$ and $\frac{v_0}{e_1}<\frac{v_1}{e_1}<\frac{2v_0}{e_1}$, $S'$ is a numerical semigroup with $\mathrm{msg}(S')=\left\{\frac{v_0}{e_1},\frac{v_1}{e_1}\right\}$, and its elements smaller than $\frac{2v_0}{e_1}$ are $0$, $\frac{v_0}{e_1}$ and $\frac{v_1}{e_1}$. Hence $\mathrm{next}_{S'}\left(\frac{v_0}{e_1}\right)=\frac{v_1}{e_1}$ and $\mathrm{next}_{S'}\left(\frac{v_1}{e_1}\right)=\frac{2v_0}{e_1}$, and Lemma~\ref{corolario concentracion generadores minimales} gives
\[
\C(S')=\frac{1}{e_1}\max\{v_1-v_0,\,2v_0-v_1\}.
\]

We now prove that $\C(S)\leq e_1\C(S')$. By definition of $\C(S')$, two consecutive elements of
$S'\setminus\{0\}$ differ by at most $\C(S')$. Multiplying by $e_1$, two consecutive elements of
$e_1S'\setminus\{0\}$, where $e_1S'=\{e_1t : t\in S'\}=\langle v_0,v_1\rangle$, differ by at most $e_1\C(S')$.

Let $s\in S\setminus\{0\}$. Since $v_0\leq s$ and $v_0\in e_1S'$, we can take the largest element $t$
of $e_1S'$ such that $t\leq s$, which is non-zero.
Let $t'$ be the element of $e_1S'$ that follows $t$. Then $t'\leq t+e_1\C(S')$, by the previous
paragraph, and $t'>s$, since otherwise $t'$ would be an element of $e_1S'$ with $t<t'\leq s$,
contradicting the choice of $t$. Moreover, $t'\in S$, because $e_1S'\subseteq S$. Thus $t'$ is
an element of $S$ greater than $s$, and, since $\mathrm{next}_S(s)$ is the smallest element
of $S$ such that $s<\mathrm{next}_S(s)$, we get
\[
\mathrm{next}_S(s)\leq t'\leq t+e_1\C(S')
\leq s+e_1\C(S').
\]
Therefore, $\mathrm{next}_S(s)-s\leq e_1\C(S')$ for every $s\in S\setminus\{0\}$, so $\C(S)\leq e_1\C(S')=\max\{v_1-v_0,\,2v_0-v_1\}$, which is~\eqref{ineq:concentration}.
\end{proof}

In Section~\ref{sec:comun} we show that \eqref{ineq:concentration} is an equality for the semigroups of values of plane branches.

\section{Concentration of the value semigroup of a plane algebroid branch}
\label{sec:comun}

Let $\Gamma(C)$ be the semigroup of values associated with a plane branch $C$ and $n=v_0<v_1<\cdots<v_g$ be its minimal system of generators (see~\cite[Chapter~II, Section~1]{Zariski-book} for $K=\mathbb{C}$ and~\cite{GB-P-2016} for $K$ of any characteristic). Set \begin{equation}\label{notacion:ek-nk}
e_k=\gcd(n,v_1,\ldots,v_k)\quad\text{and}\quad n_k=\frac{e_{k-1}}{e_k}\quad\text{for } 1\leq k\leq g,\quad\text{with } e_0=n.
\end{equation}

\noindent A semigroup $\langle v_0,v_1,\ldots,v_g\rangle$ generated by positive integers $v_0<v_1<\cdots<v_g$, with $e_k$ and $n_k$ defined by~\eqref{notacion:ek-nk} for $n=v_0$, is \emph{strongly increasing} (SI for short) if $e_0>e_1>\cdots>e_g=1$ and

\begin{equation}\label{eq:condicion SI}
n_kv_k<v_{k+1} \quad\text{ for } 1\leq k\leq g-1.
\end{equation}

\noindent The semigroup $\Gamma(C)$ is SI (see~\cite{Angermuller} and~\cite{GB-P-2016}), and, conversely, every SI-semigroup is the semigroup of
values of a plane algebroid branch (see~\cite[Theorem~6.5]{GB-P-2016}).

\begin{corollary}\label{cor:concentracion SI}
Let $\Gamma=\langle v_0,v_1,\ldots,v_g\rangle$ be an SI-semigroup with $g\geq 1$. Then
\begin{equation}\label{eq:concentracion SI}
\C(\Gamma)=\begin{cases}
\max\{v_1-v_0,\;2v_0-v_1\} & \text{if } v_1<2v_0,\\
v_0 & \text{if } v_1>2v_0.
\end{cases}
\end{equation}
\end{corollary}

\begin{proof}
If $v_1>2v_0$, then $\C(\Gamma)=v_0$ by
Proposition~\ref{proposicion la multiplicidad acota la concentracion}.

Let $v_1<2v_0$. If $g\geq 2$, then, since $\Gamma$ is strongly increasing, \eqref{eq:condicion SI} gives $v_k\geq v_2>n_1v_1\geq 2v_1>2v_0$ for every $k\geq 2$. As $v_0+v_1$ and $2v_1$ are also greater than $2v_0$, the only elements of $\Gamma$ in the interval $(0,2v_0]$ are $v_0$, $v_1$ and $2v_0$, also when $g=1$. Hence $\mathrm{next}_\Gamma(v_0)=v_1$ and $\mathrm{next}_\Gamma(v_1)=2v_0$, so $\C(\Gamma)\geq\max\{v_1-v_0,\,2v_0-v_1\}$, and equality holds by Proposition~\ref{prop:cota concentracion}.
\end{proof}

\medskip

If $v_1=hn+r$ is the Euclidean division of $v_1$ by $n$, then, by~\cite[Section~5]{Zariski},
\[
n=\mult(C^{(0)})=\mult(C^{(1)})=\cdots=\mult(C^{(h-1)}),
\]
and $\mult(C^{(h)})=r$ if $r\neq0$. Since $n<v_1$, we have $h\geq1$, and hence
\begin{equation}\label{multstrictransform}
\mult(\widetilde C)=
\begin{cases}
v_1-n & \text{if } h=1,\\
n & \text{if } h\geq2.
\end{cases}
\end{equation}

\begin{proposition}\label{proposicion concentracion curva}
Let $C$ be a singular plane algebroid branch and $n< v_1<\cdots < v_g$ the minimal system of generators of the semigroup of values $\Gamma(C)$. The concentration of $\Gamma(C)$ is
\[
\C(\Gamma(C))=\max\{\mult(\widetilde C), \mult(C)-\mult(\widetilde C)\},
\]
where $\widetilde C$ is the strict transform of $C$.
\end{proposition}

\begin{proof} 
Let $v_1=hn+r$ be the Euclidean division of $v_1$ by $n$, so that $h=1$ if $v_1<2n$ and $h\geq2$ if $v_1>2n$. By \eqref{eq:concentracion SI} and \eqref{multstrictransform}, we get
\begin{equation*}
\C(\Gamma(C))=
\begin{cases}
\max\{v_1-n,\,2n-v_1\} & \text{if } h=1\\
n & \text{if } h\geq2
\end{cases}
\;=\;\max\{\mult(\widetilde C),\,\mult(C)-\mult(\widetilde C)\}.
\end{equation*}
\end{proof}

\begin{corollary}
    The concentration of the semigroup of values of a plane algebroid branch $C$ is an invariant of the equisingularity class of $C$.
\end{corollary}

In particular the sequence $\mathcal S:=(\C(\Gamma(C^{(i)})))_{i=0}^{r-1}$ is an invariant of the equisingularity class of $C$, where $r$ is the number of steps in the resolution of the singularity of $C$, that is, $C^{(r)}$ is the first smooth strict transform of $C$.

\begin{remark}
The concentration sequence is not a complete equisingularity invariant: it is determined by the multiplicity sequence, but it does not determine it. For instance, the semigroups $\langle 5,7\rangle$ and $\langle 5,8\rangle$ have multiplicity sequences $(5,2,2,1)$ and $(5,3,2,1)$, respectively, and the same concentration sequence $\mathcal{S}=(3,2,1)$. This leads to the problem of describing all the multiplicity sequences with a given concentration sequence, which we address in Section~\ref{sec:estructura}.
\end{remark}

\begin{lemma}\label{lema cota multiplicidad y concentracion curva}
Let $C$ be a singular plane algebroid branch. Then,
\[
\frac{\mult(C)}{2}\leq\mathcal C(\Gamma(C))\leq \mult (C).
\]
\end{lemma}

\begin{proof}
The inequality $\mathcal C(\Gamma(C))\leq\mult(C)$ has already been established in Proposition~\ref{proposicion la multiplicidad acota la concentracion}.
Let $\widetilde C$ be the strict transform of $C$.
By Proposition~\ref{proposicion concentracion curva}, $\mathcal C(\Gamma(C))=\max\{\mult(\widetilde C),\
\mult(C)-\mult(\widetilde C)\}$, and the maximum of two numbers is at least their mean, which is $\frac{\mult(C)}{2}$.
\end{proof}

Any non-increasing sequence of positive integers can be grouped into maximal runs of constant value, which we call \emph{blocks}.
Concretely, given a sequence
$(m_0,m_1,\ldots,m_r)$,
its blocks are
\begin{equation}\label{eq:bloques}
(\lambda_0,\ell_0),\,
(\lambda_1,\ell_1),\,\ldots,\,
(\lambda_s,\ell_s),
\quad
\lambda_0>\lambda_1>\cdots>\lambda_s,
\end{equation}
where $\lambda_j$ denotes the common value and $\ell_j$ the length of the $j$-th run.
For example, the sequence $(4,2,2,2,2,2,1)$ gives blocks $(4,1),(2,5),(1,1)$.

We recall a result that will be used throughout this section.
Let $C$ be a plane algebroid branch
with value semigroup
$\Gamma(C)=\langle v_0,v_1,\ldots,v_g\rangle$,
where $v_0<v_1<\cdots<v_g$ is the minimal system of generators.
Let $e_i$ and $n_i$ be as in
\eqref{notacion:ek-nk}, so that
$e_0=v_0$ and $e_g=1$. Consider the Euclidean algorithms

\begin{equation}\label{eq:algoritmos}
\begin{array}{c|c}
\text{Algorithm }1 & \text{Algorithm }i\ (i=2,\ldots,g)\\[4pt]
\hline\\[-6pt]
\begin{aligned}
v_1 &= h_{1,0}\,v_0 + r_{1,1}\\
v_0 &= h_{1,1}\,r_{1,1}+r_{1,2}\\
&\vdots\\
r_{1,k_1-1}&=h_{1,k_1}\,e_1
\end{aligned}
&
\begin{aligned}
v_i-n_{i-1}v_{i-1} &= h_{i,0}\,e_{i-1}+r_{i,1}\\
e_{i-1}&=h_{i,1}\,r_{i,1}+r_{i,2}\\
&\vdots\\
r_{i,k_i-1}&=h_{i,k_i}\,e_i
\end{aligned}
\end{array}
\end{equation}

The following proposition
describes the multiplicity sequence of $C$
in terms of these algorithms.

\begin{proposition}[\cite{EnriquesChisini}, \cite{Angermuller}, {\cite[Theorem~2.8]{Zhang}}]\label{prop:Angermuller}
The multiplicity sequence of $C$
is obtained by concatenating,
for $i=1,\ldots,g$,
the pairs $(\text{value},\text{length})$
produced by the Euclidean algorithms~\eqref{eq:algoritmos}:
\begin{equation}\label{eq:pares}
\underbrace{(e_0,h_{1,0}),(r_{1,1},h_{1,1}),\ldots,(e_1,h_{1,k_1})}_{\text{algorithm }1},
\ \ldots,\
\underbrace{(e_{g-1},h_{g,0}),(r_{g,1},h_{g,1}),\ldots,(1,h_{g,k_g})}_{\text{algorithm }g},
\end{equation}
and adding the lengths of any two consecutive pairs sharing the same value.
Concretely, the last division of algorithm $i-1$
and the first division of algorithm $i$
both produce a pair of value $e_{i-1}$,
so they merge into one block of value $e_{i-1}$
and length $h_{i-1,k_{i-1}}+h_{i,0}$:
\begin{equation}\label{eq:pares-fusionados}
(e_0,h_{1,0}),(r_{1,1},h_{1,1}),\ldots,
(e_1,\,h_{1,k_1}+h_{2,0}),
(r_{2,1},h_{2,1}),\ldots,
(e_2,\,h_{2,k_2}+h_{3,0}),\ldots
\end{equation}
In~\eqref{eq:pares-fusionados} all values are strictly decreasing:
no two consecutive pairs share the same value,
so these are precisely the blocks of the multiplicity sequence.
\end{proposition}

The sequence in
Proposition~\ref{prop:Angermuller}
is the multiplicity sequence of the
minimal embedded resolution of $C$.
The multiplicity sequence considered
in this paper stops at the first smooth
strict transform $C^{(r)}$,
so it is obtained by replacing
the last pair $(1,h_{g,k_g})$ by $(1,1)$.
Since $h_{g,k_g}=r_{g,k_g-1}$
is the value of the preceding block,
no information is lost.

\begin{proposition}\label{recovery-raw}
Let $(\lambda_0,\ell_0),\ldots,(\lambda_s,\ell_s)$
be the blocks of the multiplicity sequence of $C$.
Then $v_0=\lambda_0$, $v_1=\ell_0\lambda_0+\lambda_1$,
and for $i=2,\ldots,g$,
letting $j_i$ be the index such that $\lambda_{j_i}=e_{i-1}$
and setting
$\kappa_{j_i}=\ell_{j_i}\lambda_{j_i}-\lambda_{j_i-1}+\lambda_{j_i+1}$,
one has
\[
v_i = n_{i-1}\,v_{i-1}+\kappa_{j_i}.
\]
\end{proposition}

\begin{proof}
By Proposition~\ref{prop:Angermuller},
the block at position $j_i$ has value $\lambda_{j_i}=e_{i-1}$
and merged length $\ell_{j_i}=h_{i-1,k_{i-1}}+h_{i,0}$.
The block immediately preceding it has value
$\lambda_{j_i-1}=r_{i-1,k_{i-1}-1}$,
the penultimate remainder of algorithm $i-1$;
the last step of that algorithm divides it exactly,
so $\lambda_{j_i-1}=h_{i-1,k_{i-1}}\cdot e_{i-1}$,
giving $h_{i-1,k_{i-1}}=\lambda_{j_i-1}/e_{i-1}$
and hence $h_{i,0}=\ell_{j_i}-\lambda_{j_i-1}/e_{i-1}$.
The first remainder of algorithm $i$ is $r_{i,1}=\lambda_{j_i+1}$.
Therefore
\[
h_{i,0}e_{i-1}+r_{i,1}
=\bigl(\ell_{j_i}-\tfrac{\lambda_{j_i-1}}{e_{i-1}}\bigr)e_{i-1}+\lambda_{j_i+1}
=\ell_{j_i}\lambda_{j_i}-\lambda_{j_i-1}+\lambda_{j_i+1}
=\kappa_{j_i}.
\]
Since algorithm $i$ in \eqref{eq:algoritmos}
gives $v_i-n_{i-1}v_{i-1}=h_{i,0}e_{i-1}+r_{i,1}$,
we obtain
\[
v_i = n_{i-1}\,v_{i-1}+\kappa_{j_i}.\qedhere
\]
\end{proof}

\begin{corollary}\label{cor:vi-cerrada}
Under the hypotheses of
Proposition~\ref{recovery-raw},
for $i=2,\ldots,g$,
\[
v_i
=
\Bigl(\prod_{t=1}^{i-1}n_t\Bigr)\,v_1
+
\sum_{q=2}^{i}
\Bigl(\prod_{t=q}^{i-1}n_t\Bigr)\,\kappa_{j_q},
\]
where the empty product $($for $q=i)$
equals $1$.
\end{corollary}
\begin{proof}
By induction on $i$.
For $i=2$:
$v_2=n_1 v_1+\kappa_{j_2}$,
which is Proposition~\ref{recovery-raw}.
Assume the formula holds for $i-1$.
Then Proposition~\ref{recovery-raw} gives
\begin{align*}
v_i
&=n_{i-1}\,v_{i-1}+\kappa_{j_i}\\
&=n_{i-1}\Bigl[
\Bigl(\prod_{t=1}^{i-2}n_t\Bigr)v_1
+\sum_{q=2}^{i-1}
\Bigl(\prod_{t=q}^{i-2}n_t\Bigr)
\kappa_{j_q}
\Bigr]+\kappa_{j_i}\\
&=\Bigl(\prod_{t=1}^{i-1}n_t\Bigr)v_1
+\sum_{q=2}^{i-1}
\Bigl(\prod_{t=q}^{i-1}n_t\Bigr)
\kappa_{j_q}
+\kappa_{j_i},
\end{align*}
which is the formula for $i$.
\end{proof}

Proposition~\ref{recovery-raw} expresses
the value semigroup of a branch
in terms of the blocks of its
multiplicity sequence.
Conversely, the following proposition
characterizes which sequences of blocks
are the block sequences of multiplicity
sequences of plane algebroid branches,
and recovers the corresponding
value semigroup. Although it follows from the classical description of multiplicity sequences by Euclidean algorithms recalled above, we have not found it stated in this form in the literature, so we include a proof. We point out that a characterization of the multiplicity sequences of plane branches is also given in \cite[Theorem 3.2]{Barucci2003}, although its formulation differs from the one presented here.

\begin{proposition}\label{prop:mult-to-semigrup}
Let $(\lambda_0,\ell_0),\ldots,(\lambda_s,\ell_s)$
be a sequence of blocks with
$\lambda_0>\lambda_1>\cdots>\lambda_s=1$,
$\ell_j\geq 1$ for all $j=0,\ldots,s-1$,
and $\ell_s=1$.
Let $J\subseteq\{1,\ldots,s-1\}$
be the set of indices $j$
such that $\lambda_j$ divides $\lambda_{j-1}$.
The sequence is the block sequence
of the multiplicity sequence
of a plane algebroid branch
if and only if the following two conditions hold:
\begin{enumerate}
\item[(1)] For every $j\in\{1,\ldots,s-1\}$
with $j\notin J$:
\[
\lambda_{j+1}=\lambda_{j-1}-\ell_j\,\lambda_j.
\]
\item[(2)] For every $j\in J$:
\begin{equation}\label{eq:cond-kappa}
\ell_j\,\lambda_j\geq\lambda_{j-1}.
\end{equation}
\end{enumerate}

When both conditions hold,
set $v_0=\lambda_0$ and $v_1=\ell_0\lambda_0+\lambda_1$.
If $J=\emptyset$,
then $\langle v_0,v_1\rangle$ is the value semigroup
of every plane algebroid branch
with this multiplicity sequence.
If $J\neq\emptyset$,
write $J=\{j_2<j_3<\cdots<j_g\}$,
put $j_1:=0$,
and define inductively, for $i=2,\ldots,g$,
\begin{equation}\label{eq:vi}
v_i
=\frac{\lambda_{j_{i-1}}}{\lambda_{j_i}}\,v_{i-1}
+\ell_{j_i}\lambda_{j_i}
-\lambda_{j_i-1}+\lambda_{j_i+1}.
\end{equation}
Then $\langle v_0,v_1,\ldots,v_g\rangle$
is a strongly increasing semigroup
that is the value semigroup of every plane algebroid branch
with this multiplicity sequence. Moreover, $\gcd(v_0,\ldots,v_{i-1})=\lambda_{j_i}$
for $i=1,\ldots,g$.
\end{proposition}

\begin{proof}
Throughout, for $j\in J$ we write
$\kappa_j=\ell_j\lambda_j-\lambda_{j-1}+\lambda_{j+1}$,
as in Proposition~\ref{recovery-raw}.
Recall that the multiplicity sequence
stops at the first smooth strict transform,
so its last block is $(1,1)$.

Suppose first that
$(\lambda_0,\ell_0),\ldots,(\lambda_s,\ell_s)$
are the blocks of the multiplicity sequence
of a plane algebroid branch $C$,
and let $\Gamma(C)=\langle v_0,\ldots,v_g\rangle$,
with $e_i$ and $n_i$ as in~\eqref{notacion:ek-nk}.
By the convention just recalled, $\ell_s=1$.
In algorithm~$i$ of~\eqref{eq:algoritmos}
put $r_{i,0}=e_{i-1}$, $r_{i,k_i}=e_i$
and $r_{i,k_i+1}=0$,
so that
\[
r_{i,t-1}=h_{i,t}\,r_{i,t}+r_{i,t+1}
\qquad (1\leq t\leq k_i),
\]
with $0<r_{i,t+1}<r_{i,t}$ for $t<k_i$.
By Proposition~\ref{prop:Angermuller},
the blocks are the pairs
$(r_{i,t},h_{i,t})$, $0\leq t\leq k_i$,
where for $1\leq i\leq g-1$ the pair
$(r_{i,k_i},h_{i,k_i})=(e_i,h_{i,k_i})$
is merged with $(r_{i+1,0},h_{i+1,0})=(e_i,h_{i+1,0})$.

Let $1\leq j\leq s-1$.
The value $\lambda_j$ is $r_{i,t}$ for some $i$
and some $1\leq t\leq k_i$,
and the preceding block has value
$\lambda_{j-1}=r_{i,t-1}$.
If $t<k_i$, the block is not merged,
$\ell_j=h_{i,t}$ and $\lambda_{j+1}=r_{i,t+1}$, so
\[
\lambda_{j-1}=\ell_j\lambda_j+\lambda_{j+1},
\qquad 0<\lambda_{j+1}<\lambda_j .
\]
Hence $\lambda_j$ does not divide $\lambda_{j-1}$,
that is, $j\notin J$,
and condition~(1) holds.
If $t=k_i$, then $\lambda_j=e_i$
with $i\leq g-1$ (because $\lambda_j\neq 1$).
Now $\lambda_{j-1}=r_{i,k_i-1}=h_{i,k_i}e_i$,
so $j\in J$, and
$\ell_j=h_{i,k_i}+h_{i+1,0}\geq h_{i,k_i}
=\lambda_{j-1}/\lambda_j$,
which is condition~(2).

We have shown that $J$ consists exactly
of the indices of the blocks of values
$e_1,\ldots,e_{g-1}$.
Thus $J=\{j_2<\cdots<j_g\}$ with
$\lambda_{j_i}=e_{i-1}$,
and also $\lambda_{j_1}=\lambda_0=e_0$.
By Proposition~\ref{recovery-raw},
$v_0=\lambda_0$, $v_1=\ell_0\lambda_0+\lambda_1$ and
$v_i=n_{i-1}v_{i-1}+\kappa_{j_i}$ for $i\geq 2$,
where
$n_{i-1}=e_{i-2}/e_{i-1}=\lambda_{j_{i-1}}/\lambda_{j_i}$.
This is~\eqref{eq:vi}.
Hence the value semigroup of every branch
with this multiplicity sequence
is the one given in the statement,
and $\gcd(v_0,\ldots,v_{i-1})=e_{i-1}=\lambda_{j_i}$.

Conversely, assume that $\ell_s=1$
and that conditions~(1) and~(2) hold.
We show that the blocks are those of
the multiplicity sequence of a branch.
Let $g=|J|+1$, $J=\{j_2<\cdots<j_g\}$,
$j_1=0$ and $j_{g+1}=s$,
and define $v_0,\ldots,v_g$ as in the statement.

(a) For $1\leq i\leq g$,
$\gcd(\lambda_{j_i},\lambda_{j_i+1})=\lambda_{j_{i+1}}$.
Indeed, for $j_i<j<j_{i+1}$ we have
$1\leq j\leq s-1$ and $j\notin J$,
so condition~(1) gives
$\lambda_{j-1}=\ell_j\lambda_j+\lambda_{j+1}$,
and hence
$\gcd(\lambda_{j-1},\lambda_j)=\gcd(\lambda_j,\lambda_{j+1})$.
Applying this for $j=j_i+1,\ldots,j_{i+1}-1$
(no step is needed if $j_{i+1}=j_i+1$) gives
$\gcd(\lambda_{j_i},\lambda_{j_i+1})
=\gcd(\lambda_{j_{i+1}-1},\lambda_{j_{i+1}})
=\lambda_{j_{i+1}}$,
because $\lambda_{j_{i+1}}$ divides
$\lambda_{j_{i+1}-1}$
(as $j_{i+1}\in J$ if $i<g$,
and $\lambda_{j_{g+1}}=\lambda_s=1$).

(b) Put $e_i=\gcd(v_0,\ldots,v_i)$.
Then $e_i=\lambda_{j_{i+1}}$ for $0\leq i\leq g$.
For $i=0$, $e_0=v_0=\lambda_0=\lambda_{j_1}$.
For $i=1$,
$e_1=\gcd(\lambda_0,\ell_0\lambda_0+\lambda_1)
=\gcd(\lambda_0,\lambda_1)=\lambda_{j_2}$ by~(a).
Let $i\geq 2$ and assume $e_{i-1}=\lambda_{j_i}$.
Since $\lambda_{j_i}$ divides $v_{i-1}$,
it divides the first summand of~\eqref{eq:vi}.
Since $j_i\in J$, $\lambda_{j_i}$ divides
$\ell_{j_i}\lambda_{j_i}-\lambda_{j_i-1}$,
so $\kappa_{j_i}\equiv\lambda_{j_i+1}
\pmod{\lambda_{j_i}}$.
Hence
$e_i=\gcd(\lambda_{j_i},v_i)
=\gcd(\lambda_{j_i},\lambda_{j_i+1})
=\lambda_{j_{i+1}}$ by~(a).
In particular
$e_0>e_1>\cdots>e_g=\lambda_s=1$,
$n_i:=e_{i-1}/e_i=\lambda_{j_i}/\lambda_{j_{i+1}}$,
and~\eqref{eq:vi} reads
$v_i=n_{i-1}v_{i-1}+\kappa_{j_i}$.

(c) By condition~(2),
$\kappa_{j_i}\geq\lambda_{j_i+1}>0$,
so $v_i>n_{i-1}v_{i-1}$ for $2\leq i\leq g$.
Together with~(b), this shows that
$\langle v_0,\ldots,v_g\rangle$
is a strongly increasing semigroup,
hence, by \cite[Theorem~6.5]{GB-P-2016},
the value semigroup of a plane algebroid branch $C$.

(d) We compute the algorithms~\eqref{eq:algoritmos}
for $\Gamma(C)=\langle v_0,\ldots,v_g\rangle$.
Algorithm~$1$ starts with
$v_1=\ell_0v_0+\lambda_1$, $0<\lambda_1<\lambda_0$,
so $h_{1,0}=\ell_0$ and $r_{1,1}=\lambda_1$;
by condition~(1) the following divisions are
$\lambda_{j-1}=\ell_j\lambda_j+\lambda_{j+1}$
for $1\leq j<j_2$,
and the last one is
$\lambda_{j_2-1}=(\lambda_{j_2-1}/\lambda_{j_2})\,\lambda_{j_2}$.
For $2\leq i\leq g$, algorithm~$i$ starts with
\[
v_i-n_{i-1}v_{i-1}=\kappa_{j_i}
=\Bigl(\ell_{j_i}-\frac{\lambda_{j_i-1}}{\lambda_{j_i}}\Bigr)
\lambda_{j_i}+\lambda_{j_i+1},
\]
where the coefficient of $\lambda_{j_i}$
is $\geq 0$ by condition~(2)
and $0<\lambda_{j_i+1}<\lambda_{j_i}$;
so $h_{i,0}=\ell_{j_i}-\lambda_{j_i-1}/\lambda_{j_i}$
and $r_{i,1}=\lambda_{j_i+1}$,
and the following divisions are given by
condition~(1) as before,
up to
$\lambda_{j_{i+1}-1}
=(\lambda_{j_{i+1}-1}/\lambda_{j_{i+1}})\,\lambda_{j_{i+1}}$.
By Proposition~\ref{prop:Angermuller},
the multiplicity sequence of $C$
has a block of value $\lambda_j$ and length $\ell_j$
for every $j\notin J\cup\{s\}$,
and, for $j=j_i\in J$,
a merged block of value $\lambda_{j_i}$ and length
$\lambda_{j_i-1}/\lambda_{j_i}+h_{i,0}=\ell_{j_i}$.
Its last block is $(1,1)=(\lambda_s,\ell_s)$.
Hence the multiplicity sequence of $C$
has blocks
$(\lambda_0,\ell_0),\ldots,(\lambda_s,\ell_s)$.
\end{proof}

\begin{example}\label{ex:mult-to-semigrup}
We illustrate Proposition~\ref{prop:mult-to-semigrup}
with four examples.

\emph{(1) Case $J=\emptyset$.}
Consider the sequence $(3,3,3,1)$
with blocks $(\lambda_0,\ell_0)=(3,3)$
and $(\lambda_1,\ell_1)=(1,1)$.
Here $s=1$, so $\{1,\ldots,s-1\}=\emptyset$ and $J=\emptyset$.
Then $v_0=3$ and $v_1=3\cdot 3+1=10$.
By Proposition~\ref{prop:mult-to-semigrup},
the semigroup $\langle 3,10\rangle$
is the value semigroup
of every plane algebroid branch
with multiplicity sequence $(3,3,3,1)$.

\emph{(2) Case $J\neq\emptyset$.}
Consider the sequence $(6,6,3,3,3,1)$
with blocks
$(\lambda_0,\ell_0)=(6,2)$,
$(\lambda_1,\ell_1)=(3,3)$,
$(\lambda_2,\ell_2)=(1,1)$.
Here $s=2$ and $\{1,\ldots,s-1\}=\{1\}$.
Since $3$ divides $6$, we have $J=\{1\}$,
so condition~(1) is vacuous.
Condition~(2) at $j=1$:
$\ell_1\lambda_1=3\cdot 3=9\geq 6=\lambda_0$,
so it holds.
We compute $v_0=6$, $v_1=2\cdot 6+3=15$,
and by~\eqref{eq:vi} with $j_1=0$
and $j_2=1$:
$v_2=(\lambda_0/\lambda_1)\,v_1
+\ell_1\lambda_1-\lambda_0+\lambda_2
=(6/3)\cdot 15+9-6+1=34$.
The value semigroup of every plane algebroid branch
with multiplicity sequence $(6,6,3,3,3,1)$
is $\langle 6,15,34\rangle$.

\emph{(3) Failure of condition~(2).}
Consider the blocks
$(\lambda_0,\ell_0)=(4,2)$,
$(\lambda_1,\ell_1)=(2,1)$,
$(\lambda_2,\ell_2)=(1,1)$.
Here $J=\{1\}$ since $2$ divides $4$,
so condition~(2) applies at $j=1$:
$\ell_1\lambda_1=1\cdot 2=2<4=\lambda_0$.
Hence this is not the block sequence
of any plane algebroid branch.

\emph{(4) Failure of condition~(1).}
Consider the blocks
$(\lambda_0,\ell_0)=(5,1)$,
$(\lambda_1,\ell_1)=(3,1)$,
$(\lambda_2,\ell_2)=(1,1)$.
Here $\{1,\ldots,s-1\}=\{1\}$
and $3$ does not divide $5$,
so $J=\emptyset$
and condition~(1) applies at $j=1$:
it requires $\lambda_2=\lambda_0-\ell_1\lambda_1=5-3=2$,
but $\lambda_2=1\neq 2$.
Hence this is not the block sequence
of any plane algebroid branch.
\end{example}

\begin{lemma}\label{lem:bloques-divide}
Let $(\lambda_0,\ell_0),\ldots,(\lambda_s,\ell_s)$
be the blocks of the multiplicity sequence
of a plane algebroid branch, and let
$1\leq j\leq s$ be such that $\lambda_j$
divides $\lambda_{j-1}$. Then $\lambda_j$
divides $\lambda_i$ for every $i\leq j$.
\end{lemma}

\begin{proof}
The integer $\lambda_j$ divides $\lambda_j$
and $\lambda_{j-1}$. If $\lambda_j$ divides
$\lambda_i$ and $\lambda_{i+1}$ for some
$1\leq i\leq j-1$, then it divides
$\lambda_{i-1}$, because either $\lambda_i$
divides $\lambda_{i-1}$, or
$\lambda_{i-1}=\ell_i\lambda_i+\lambda_{i+1}$
by condition~(1) of
Proposition~\ref{prop:mult-to-semigrup}.
The lemma follows by descending induction
on $i$.
\end{proof}

In terms of
Proposition~\ref{prop:Angermuller}, the
values $\lambda_j$ with $1\leq j\leq s$ that
divide $\lambda_{j-1}$ are exactly
$e_1,\ldots,e_g$, the last nonzero remainders
of the Euclidean algorithms~\eqref{eq:algoritmos};
thus Lemma~\ref{lem:bloques-divide} expresses
the fact that $e_i$ divides all the
multiplicities preceding the block of
value $e_i$.

\begin{lemma}\label{lem:delante}
Let $(m_0,\ldots,m_{r-1},1)$
be the multiplicity sequence
of a plane algebroid branch, with $r\geq1$,
let $\ell$ be the number of initial copies
of $m_0$ in it, and let $\rho$ be the term
that follows them.
For an integer $m'\geq m_0$,
the sequence $(m',m_0,\ldots,m_{r-1},1)$
is the multiplicity sequence
of a plane algebroid branch
if and only if
\[
m'\in\{m_0,\,2m_0,\,\ldots,\,\ell m_0,\ \ell m_0+\rho\}.
\]
\end{lemma}

\begin{proof}
Since $r\geq1$, we have $m_0\geq m_{r-1}\geq2$,
so $\rho$ exists and $0<\rho<m_0$.
Let $(\lambda_0,\ell_0),\ldots,(\lambda_s,\ell_s)$
be the blocks of $(m_0,\ldots,m_{r-1},1)$,
so that $(\lambda_0,\ell_0)=(m_0,\ell)$,
$\lambda_1=\rho$ and $s\geq1$.
By Proposition~\ref{prop:mult-to-semigrup},
these blocks satisfy conditions~(1) and~(2),
and $(m',m_0,\ldots,m_{r-1},1)$
is a multiplicity sequence
if and only if its own blocks satisfy them.
If $m'=m_0$, its blocks are
$(m_0,\ell+1),(\lambda_1,\ell_1),\ldots,(\lambda_s,\ell_s)$;
since the length of the first block appears
in no condition, the conditions are the old
ones, which hold.
If $m'>m_0$, its blocks are
$(m',1),(\lambda_0,\ell_0),\ldots,(\lambda_s,\ell_s)$,
so every old block moves one position to the
right. Each condition only involves a block
and its two neighbours, so the only new
condition is the one at the block $(m_0,\ell)$,
which is now in position~$1$ and has
neighbours $m'$ and $\rho$:
if $m_0$ divides $m'$, it is condition~(2),
that is, $m'\leq\ell m_0$;
if $m_0$ does not divide $m'$, it is
condition~(1), that is, $m'=\ell m_0+\rho$.
Since $0<\rho<m_0$, the number $\ell m_0+\rho$
is not a multiple of $m_0$, and the lemma
follows.
\end{proof}

\begin{lemma}\label{lem:tipoB}
Let $m=(m_0,\ldots,m_{r-1},1)$
be the multiplicity sequence
of a plane algebroid branch $C$,
with concentration sequence
$\mathcal{S}=(c_0,\ldots,c_{r-1})$,
and let $0\leq i\leq r-2$.
\begin{enumerate}
\item[(a)] If $c_i=m_{i+1}$, then
$m_i\in\{\,m_{i+1},\ m_{i+1}+m_{i+2}\,\}$.
\item[(b)] If $c_i>m_{i+1}$, then
$m_{i+2}=m_{i+1}$.
\end{enumerate}
\end{lemma}

\begin{proof}
Let $(\lambda_0,\ell_0),\ldots,(\lambda_s,\ell_s)$
be the blocks of $m$.
By Proposition~\ref{proposicion concentracion curva}
applied to $C^{(i)}$,
$c_i=\max\{m_{i+1},m_i-m_{i+1}\}$.
Since $C^{(r)}$ is the first smooth
strict transform, $m_{r-1}\geq 2$,
so $m_{i+1}\geq 2$ because $i+1\leq r-1$.
If $m_i>m_{i+1}$, the block of value
$m_{i+1}$ starts at position $i+1$;
let $j$ be its index and $\ell$ its length.
Then $1\leq j\leq s-1$,
because this block is preceded
by the block of value $m_i$
and $m_{i+1}>1=\lambda_s$,
and by Proposition~\ref{prop:mult-to-semigrup}
conditions~(1) and~(2) hold for these blocks.

(a) Since $c_i=m_{i+1}$, we get
$m_i-m_{i+1}\leq m_{i+1}$, that is,
$m_{i+1}\leq m_i\leq 2m_{i+1}$.
Assume $m_i\neq m_{i+1}$.
If $m_{i+1}$ does not divide $m_i$,
condition~(1) gives
$\lambda_{j+1}=m_i-\ell m_{i+1}>0$;
since $m_i\leq 2m_{i+1}$ this forces $\ell=1$,
so $m_{i+2}=\lambda_{j+1}=m_i-m_{i+1}$.
If $m_{i+1}$ divides $m_i$,
then $m_i=2m_{i+1}$
and condition~(2) gives
$\ell m_{i+1}\geq 2m_{i+1}$,
so $\ell\geq 2$ and
$m_{i+2}=m_{i+1}=m_i-m_{i+1}$.
In both cases $m_i=m_{i+1}+m_{i+2}$.

(b) Since $c_i>m_{i+1}$, we get
$c_i=m_i-m_{i+1}$, so $m_i>2m_{i+1}$.
If $m_{i+1}$ divides $m_i$,
condition~(2) gives
$\ell m_{i+1}\geq m_i>2m_{i+1}$,
so $\ell\geq 3$.
Otherwise condition~(1) gives
$0<m_i-\ell m_{i+1}<m_{i+1}$,
so $\ell m_{i+1}>m_i-m_{i+1}>m_{i+1}$
and $\ell\geq 2$.
In both cases $m_{i+2}=m_{i+1}$.
\end{proof}

\section{Multiplicity sequences with a given concentration sequence}
\label{sec:estructura}

Let $\mathcal{S}=(c_0,\ldots,c_{r-1})$ be a
sequence of positive integers and put $c_r=1$. We denote by
$\mathcal{M}(\mathcal{S})$ the set of
multiplicity sequences of plane algebroid
branches with concentration sequence
$\mathcal{S}$, and by $N(\mathcal{S})$ its
number of elements.
In this section we describe
$\mathcal{M}(\mathcal{S})$.
Theorem~\ref{teo:descripcion}
characterizes its elements by local
conditions;
Lemmas~\ref{lem:dif-aislada},
\ref{lem:tramo-inicial}
and~\ref{lem:divisibilidad} describe how
two of its elements can differ; and
Theorem~\ref{teo:producto} shows that
$\mathcal{M}(\mathcal{S})$ is a product and
computes $N(\mathcal{S})$. Finally, Proposition~\ref{prop:imagen} characterizes the sequences of positive integers that are concentration sequences of plane branches, reading them from right to left.

Throughout, the \emph{position} $k$ of a
sequence $m=(m_0,\ldots,m_{r-1},1)$ is the
place of the term $m_k$; thus positions are
numbered from $0$, and the position $k$ is
the $(k+1)$-th term of $m$.

The following theorem describes exactly
the multiplicity sequences with a given
concentration sequence. For a sequence
$m=(m_0,\ldots,m_{r-1},1)$ and
$0\leq i\leq r-1$, we denote by $\ell_i(m)$
the number of consecutive terms of $m$ equal
to $m_i$, starting with the term $m_i$, and
by $\rho_i(m)$ the first term after them;
that is, $\ell_i(m)$ and $\rho_i(m)$ are the
numbers $\ell$ and $\rho$ of
Lemma~\ref{lem:delante} for the sequence
$(m_i,\ldots,m_{r-1},1)$.
For instance, for $m=(7,7,3,3,1)$ we have
$\ell_0(m)=2$ and $\rho_0(m)=3$, while
$\ell_1(m)=1$ and $\rho_1(m)=3$.
If $i=0$ or $m_{i-1}>m_i$, then a block
$(\lambda_j,\ell_j)$ of $m$ starts at the
position $i$, and $\ell_i(m)=\ell_j$ and
$\rho_i(m)=\lambda_{j+1}$; this is the case
in~(A) below for the position $i+1$, since
$m_i>m_{i+1}$.

\begin{theorem}\label{teo:descripcion}
Let $\mathcal{S}=(c_0,\ldots,c_{r-1})$ be a
sequence of positive integers. A sequence
$m=(m_0,\ldots,m_{r-1},1)$ is the
multiplicity sequence of a plane algebroid
branch with concentration sequence
$\mathcal{S}$ if and only if
$m_{r-1}=c_{r-1}+1$ and, for every
$0\leq i\leq r-2$, writing $\ell=\ell_{i+1}(m)$
and $\rho=\rho_{i+1}(m)$, one of the following
holds:
\begin{enumerate}
\item[(A)] $c_i>m_{i+1}$, $m_i=c_i+m_{i+1}$,
and either $m_{i+1}$ divides $c_i$ and
$c_i\leq(\ell-1)m_{i+1}$, or
$c_i=(\ell-1)m_{i+1}+\rho$;
\item[(B)] $c_i=m_{i+1}$ and
$m_i\in\{m_{i+1},\ m_{i+1}+m_{i+2}\}$.
\end{enumerate}
\end{theorem}

\begin{proof}
Let $m$ be the multiplicity sequence of a
plane algebroid branch $C$ with concentration
sequence $\mathcal{S}$. Since $C^{(r)}$ is the
first smooth strict transform, $m_{r-1}\geq 2$,
and Proposition~\ref{proposicion concentracion curva}
applied to $C^{(r-1)}$ gives
$c_{r-1}=\max\{1,m_{r-1}-1\}=m_{r-1}-1$.
Let $0\leq i\leq r-2$. The sequence
$(m_{i+1},\ldots,m_{r-1},1)$ is the
multiplicity sequence of $C^{(i+1)}$, and
$(m_i,m_{i+1},\ldots,m_{r-1},1)$ is that of
$C^{(i)}$. Hence
Lemma~\ref{lem:delante}, applied to the
multiplicity sequence of $C^{(i+1)}$ and to
$m'=m_i\geq m_{i+1}$, gives
$m_i\in\{jm_{i+1}:1\leq j\leq\ell\}
\cup\{\ell m_{i+1}+\rho\}$, where
$\ell=\ell_{i+1}(m)$ and
$\rho=\rho_{i+1}(m)$ are the numbers $\ell$
and $\rho$ of
Lemma~\ref{lem:delante} for this
sequence.
By Proposition~\ref{proposicion concentracion curva},
$c_i=\max\{m_{i+1},m_i-m_{i+1}\}$; hence
$c_i=m_{i+1}$ if $m_i=m_{i+1}$ or $m_i=2m_{i+1}$,
$c_i=(j-1)m_{i+1}>m_{i+1}$ if $m_i=jm_{i+1}$
with $3\leq j\leq\ell$, and, if
$m_i=\ell m_{i+1}+\rho$, then $c_i=m_{i+1}$
when $\ell=1$ and $c_i=(\ell-1)m_{i+1}+\rho>m_{i+1}$
when $\ell\geq 2$.
If $c_i=m_{i+1}$, then either
$m_i=m_{i+1}$, or $m_i=2m_{i+1}$ with
$\ell\geq 2$, or $m_i=m_{i+1}+\rho$ with
$\ell=1$. If $\ell\geq 2$, then at least two
consecutive terms are equal to $m_{i+1}$,
that is, $m_{i+2}=m_{i+1}$, so
$2m_{i+1}=m_{i+1}+m_{i+2}$. If $\ell=1$,
then the first term after $m_{i+1}$ is
$m_{i+2}$, that is, $\rho=m_{i+2}$, so
$m_{i+1}+\rho=m_{i+1}+m_{i+2}$.
In all cases
$m_i\in\{m_{i+1},m_{i+1}+m_{i+2}\}$,
which is~(B).
If $c_i>m_{i+1}$, then
Proposition~\ref{proposicion concentracion curva}
gives $c_i=m_i-m_{i+1}$, that is,
$m_i=c_i+m_{i+1}$. Moreover, by the
computation of $c_i$ above, $c_i>m_{i+1}$
only happens in two cases: either
$m_i=jm_{i+1}$ with $3\leq j\leq\ell$, and
then $c_i=(j-1)m_{i+1}$, so $m_{i+1}$
divides $c_i$ and
$c_i\leq(\ell-1)m_{i+1}$; or
$m_i=\ell m_{i+1}+\rho$ with $\ell\geq 2$,
and then $c_i=(\ell-1)m_{i+1}+\rho$.
This is~(A).

Conversely, let $m$ satisfy the conditions
of the statement. We show by descending
induction on $i$ that $(m_i,\ldots,m_{r-1},1)$
is the multiplicity sequence of a plane
algebroid branch with concentration sequence
$(c_i,\ldots,c_{r-1})$.
For $i=r-1$, the sequence $(c_{r-1}+1,1)$
has blocks $(c_{r-1}+1,1),(1,1)$, which satisfy
the conditions of
Proposition~\ref{prop:mult-to-semigrup}
trivially, and
Proposition~\ref{proposicion concentracion curva}
gives the concentration
$\max\{1,c_{r-1}\}=c_{r-1}$.
Let $0\leq i\leq r-2$ and assume the claim for
$i+1$. Then $(m_{i+1},\ldots,m_{r-1},1)$ is
the multiplicity sequence of a plane
algebroid branch; in particular it is
non-increasing, so $m_{i+2}\leq m_{i+1}$.
Write $\ell=\ell_{i+1}(m)$ and
$\rho=\rho_{i+1}(m)$. We first show that
$m_i\geq m_{i+1}$, that $m_i$ belongs to the
set
$\{jm_{i+1}:1\leq j\leq\ell\}\cup\{\ell m_{i+1}+\rho\}$
of Lemma~\ref{lem:delante}, and that
$\max\{m_{i+1},m_i-m_{i+1}\}=c_i$.

In case~(B),
$m_i\in\{m_{i+1},m_{i+1}+m_{i+2}\}$, so
$m_i\geq m_{i+1}$. If $m_i=m_{i+1}$, then
$m_i$ is the element of the set with $j=1$.
If $m_i=m_{i+1}+m_{i+2}$, then, as in the
first part of the proof, $m_{i+2}=m_{i+1}$
when $\ell\geq 2$, and $m_{i+2}=\rho$ when
$\ell=1$; hence either $m_i=2m_{i+1}$ with
$\ell\geq 2$, or $m_i=\ell m_{i+1}+\rho$ with
$\ell=1$, and in both cases $m_i$ belongs to
the set. Moreover,
$m_i-m_{i+1}\in\{0,m_{i+2}\}$ and
$m_{i+2}\leq m_{i+1}$, so
$\max\{m_{i+1},m_i-m_{i+1}\}=m_{i+1}=c_i$.

In case~(A), $m_i=c_i+m_{i+1}>m_{i+1}$.
If $m_{i+1}$ divides $c_i$ and
$c_i\leq(\ell-1)m_{i+1}$, write
$c_i=(j-1)m_{i+1}$; since $c_i>m_{i+1}$, we
have $2\leq j-1\leq\ell-1$, that is,
$3\leq j\leq\ell$, and
$m_i=c_i+m_{i+1}=jm_{i+1}$.
If $c_i=(\ell-1)m_{i+1}+\rho$, then
$m_i=c_i+m_{i+1}=\ell m_{i+1}+\rho$.
In both cases $m_i$ belongs to the set, and
$\max\{m_{i+1},m_i-m_{i+1}\}
=\max\{m_{i+1},c_i\}=c_i$.

Now, by the induction hypothesis,
$(m_{i+1},\ldots,m_{r-1},1)$ is the
multiplicity sequence of a plane algebroid
branch $C'$ with concentration sequence
$(c_{i+1},\ldots,c_{r-1})$.
By Lemma~\ref{lem:delante}, applied to
this sequence and to $m'=m_i\geq m_{i+1}$,
the sequence $(m_i,\ldots,m_{r-1},1)$ is the
multiplicity sequence of a plane algebroid
branch $C$.
The first strict transform of $C$ has
multiplicity sequence
$(m_{i+1},\ldots,m_{r-1},1)$, so, by
Proposition~\ref{proposicion concentracion curva},
it has the same concentration sequence as
$C'$, namely $(c_{i+1},\ldots,c_{r-1})$.
Finally,
Proposition~\ref{proposicion concentracion curva}
gives
$\C(\Gamma(C))=\max\{m_{i+1},m_i-m_{i+1}\}=c_i$.
Hence the concentration sequence of $C$ is
$(c_i,\ldots,c_{r-1})$.
\end{proof}

Note that no global condition on the
blocks of $m$ has to be checked: in case~(B) both values of $m_i$
always give a multiplicity sequence, and in
case~(A) the extra condition only involves
$m_{i+1}$, $\ell_{i+1}(m)$ and $\rho_{i+1}(m)$.
Theorem~\ref{teo:descripcion} therefore gives
the following algorithm, which computes all the
multiplicity sequences with concentration
sequence $\mathcal{S}=(c_0,\ldots,c_{r-1})$.

\begin{corollary}\label{cor:algoritmo}
Start with the sequence $(c_{r-1}+1,1)$.
For $i=r-2,\ldots,0$, replace each sequence
$(m_{i+1},\ldots,m_{r-1},1)$ obtained so far
by the sequences $(m_i,m_{i+1},\ldots,m_{r-1},1)$
given by Theorem~\ref{teo:descripcion}:
two sequences if $c_i=m_{i+1}$, at most one if
$c_i>m_{i+1}$, and none if $c_i<m_{i+1}$.
The sequences obtained for $i=0$ are exactly
the multiplicity sequences of plane algebroid
branches with concentration sequence
$\mathcal{S}$.
\end{corollary}

\begin{remark}\label{rem:algoritmo}
Let $\mathcal{S}=(c_0,\ldots,c_{r-1})$ be a
sequence of positive integers, let
$m\in\mathcal{M}(\mathcal{S})$ be the
multiplicity sequence of a plane algebroid
branch $C$, and let $0\leq i\leq r-2$.
Since $m_{i+1}\geq m_{r-1}\geq 2$, the strict
transform $C^{(i+1)}$ is singular, and
Lemma~\ref{lema cota multiplicidad y concentracion curva}
applied to it gives
$m_{i+1}/2\leq c_{i+1}\leq m_{i+1}$.
Moreover,
Proposition~\ref{proposicion concentracion curva}
gives $c_i=\max\{m_{i+1},m_i-m_{i+1}\}
\geq m_{i+1}$. Hence:
\begin{enumerate}
\item[(a)] $c_{i+1}\leq m_{i+1}\leq c_i$. In
particular, $\mathcal{S}$ is non-increasing,
and $\mathcal{M}(\mathcal{S})=\emptyset$
whenever $\mathcal{S}$ is not
non-increasing.
\item[(b)] If case~(B) of
Theorem~\ref{teo:descripcion} holds at
$i$, then $c_i=m_{i+1}\leq 2c_{i+1}$. Hence,
if $c_i>2c_{i+1}$, then case~(A) holds at $i$
for every $m\in\mathcal{M}(\mathcal{S})$.
\end{enumerate}
When implementing the algorithm of
Corollary~\ref{cor:algoritmo}, (a)
allows one to discard at once every
$\mathcal{S}$ that is not non-increasing,
and~(b) shows that, at each $i$ with
$c_i>2c_{i+1}$, only case~(A) has to be
considered.
\end{remark}

\begin{example}\label{ex:arbol}
Let $\mathcal{S}=(7,4,3,2)$. We start with
$(3,1)$. For $i=2$, $c_2=3=m_3$, so we get
$(3,3,1)$ and $(4,3,1)$.
For $i=1$ and $(3,3,1)$: $c_1=4>3$, $\ell=2$,
$\rho=1$, and $4=(\ell-1)\cdot 3+\rho$, so we get
$(7,3,3,1)$. For $(4,3,1)$: $c_1=4=m_2$, so we
get $(4,4,3,1)$ and $(7,4,3,1)$.
For $i=0$ and $(7,3,3,1)$: $c_0=7=m_1$, so we
get $(7,7,3,3,1)$ and $(10,7,3,3,1)$. For
$(4,4,3,1)$: $c_0=7>4$, $\ell=2$, $\rho=3$, and
$7=4+3$, so we get $(11,4,4,3,1)$. For
$(7,4,3,1)$: $c_0=7=m_1$, so we get
$(7,7,4,3,1)$ and $(11,7,4,3,1)$.
Hence $\mathcal{M}(\mathcal{S})$ consists
of these five sequences.
By Proposition~\ref{prop:mult-to-semigrup},
their value semigroups are
$\langle 7,17\rangle$, $\langle 10,17\rangle$,
$\langle 11,15\rangle$, $\langle 7,18\rangle$
and $\langle 11,18\rangle$, respectively.
Since these semigroups are distinct, the
five branches are not equisingular: the
concentration sequence alone does not
determine the equisingularity class of a
branch. We come back to this in
Section~\ref{sec:extendida}.

Not every choice in the algorithm leads to
a multiplicity sequence.
For $\mathcal{S}=(5,2,2,1)$ we start with
$(2,1)$. For $i=2$, $c_2=2=m_3$, so we get
$(2,2,1)$ and $(3,2,1)$.
For $i=1$ and $(2,2,1)$: $c_1=2=m_2$, so we
get $(2,2,2,1)$ and $(4,2,2,1)$.
For $(3,2,1)$: $c_1=2<3=m_2$, so this
sequence is discarded.
For $i=0$ and $(2,2,2,1)$: $c_0=5>2$,
$\ell=3$, $\rho=1$, and
$5=(\ell-1)\cdot 2+\rho$, so we get
$(7,2,2,2,1)$.
For $(4,2,2,1)$: $c_0=5>4$, $\ell=1$ and
$\rho=2$; here $4$ does not divide $5$ and
$5\neq(\ell-1)\cdot 4+\rho=2$, so this
sequence is discarded: indeed,
$(9,4,2,2,1)$ does not satisfy
condition~(1) of
Proposition~\ref{prop:mult-to-semigrup}.
Hence $(7,2,2,2,1)$ is the only multiplicity
sequence with $\mathcal{S}=(5,2,2,1)$.
\end{example}

The following lemmas describe how two
multiplicity sequences with the same
concentration sequence can differ: apart from
an initial segment of length at most three,
they can only differ at isolated positions,
each of a very rigid form
(Lemmas~\ref{lem:dif-aislada}
and~\ref{lem:tramo-inicial}), and the values
before each such position are multiples of a
common value, which is one of the
integers~$e_j$
(Lemma~\ref{lem:divisibilidad}).

Let $m=(m_0,\ldots,m_{r-1},1)$ be the
multiplicity sequence of a plane algebroid
branch with concentration sequence
$\mathcal{S}=(c_0,\ldots,c_{r-1})$.
For $0\leq i\leq r-2$ we say that $i$ is
\emph{of type~A} (respectively \emph{B})
for $m$ if case~(A) (respectively~(B))
of Theorem~\ref{teo:descripcion}
holds at $i$.
Since
$c_i=\max\{m_{i+1},m_i-m_{i+1}\}\geq m_{i+1}$
by Proposition~\ref{proposicion concentracion curva},
exactly one of the two types occurs at
each $i$: $i$ is of type~A for $m$ if and
only if $m_i>2m_{i+1}$, and of type~B
if and only if $m_i\leq 2m_{i+1}$.

\begin{lemma}\label{lem:dif-aislada}
Let $m=(m_0,\ldots,m_{r-1},1)$ and
$\tilde m=(\tilde m_0,\ldots,\tilde m_{r-1},1)$
be the multiplicity sequences of two plane
algebroid branches with the same
concentration sequence
$\mathcal{S}=(c_0,\ldots,c_{r-1})$, and let
$t$ be the least index with $m_t=\tilde m_t$.
\begin{enumerate}
\item[(a)] Suppose that
$m_{k-1}=\tilde m_{k-1}=:\lambda$ and
$p:=m_k>q:=\tilde m_k$ for some $k\geq 1$.
Then
\[
\lambda=p+q,\quad
m_{k+1}=\tilde m_{k+1}=q,\quad
d:=p-q=m_{k+2}\leq q.
\]
\item[(b)] After the position $t$, the
sequences $m$ and $\tilde m$ can only differ
at isolated positions: if they differ at a
position $k>t$, then they agree at $k-1$
and at $k+1$, and the common value
$\lambda=m_{k-1}=\tilde m_{k-1}$ satisfies
$\lambda=p+q$ with the notation of~(a).
\end{enumerate}
\end{lemma}

\begin{proof}
(a) By Proposition~\ref{proposicion concentracion curva}
applied in $m$,
$c_{k-1}=\max\{p,\lambda-p\}$, and
$c_{k-1}=\max\{q,\lambda-q\}$ applied in $\tilde m$.
Since $p\neq q$, the two maxima cannot both
be attained at the first entry, nor both
at the second one; hence $p=\lambda-q$
or $\lambda-p=q$, that is, $\lambda=p+q$
in both cases, and $q<\lambda/2<p$
because $p>q$.
By Theorem~\ref{teo:descripcion},
$m_{r-1}=\tilde m_{r-1}=c_{r-1}+1$, and
$m_r=\tilde m_r=1$; since $m$ and $\tilde m$
differ at $k$, we get $k\leq r-2$.
As $\tilde m$ is non-increasing and
$\tilde m_{r-1}\geq 2$, it follows that
$q=\tilde m_k\geq\tilde m_{r-1}\geq 2$.
Let $(\lambda_0,\ell_0),\ldots,(\lambda_s,\ell_s)$
be the blocks of $m$, so that $\lambda_s=1$.
The block of value $p$ starts
at $k$, because $m_{k-1}=\lambda>p$;
let $(\lambda_j,\ell_j)$ be this block,
so that $\lambda_j=p$ and $\lambda_{j-1}=\lambda$.
Then $j\geq 1$, because the block
containing position $k-1$ has value
$\lambda\neq p$ and precedes this one,
and $j\leq s-1$ because
$p>q\geq 2>1=\lambda_s$.
Moreover, $p$ does not divide $\lambda$
because $p<\lambda<2p$.
Hence condition~(1) of
Proposition~\ref{prop:mult-to-semigrup}
gives $\lambda=\ell_j\,p+\lambda_{j+1}$
with $0<\lambda_{j+1}<p$, so $\ell_j$ and
$\lambda_{j+1}$ are the quotient and the
remainder of the division of $\lambda$ by $p$.
Since $\lambda=p+q$ with $0<q<p$, we get
$\ell_j=1$ and
$m_{k+1}=\lambda_{j+1}=q$.
By
Proposition~\ref{proposicion concentracion curva}
applied in $m$,
$c_k=\max\{m_{k+1},m_k-m_{k+1}\}
=\max\{q,p-q\}\geq q$, while, applied
in $\tilde m$,
$c_k=\max\{\tilde m_{k+1},\tilde m_k-\tilde m_{k+1}\}
=\max\{\tilde m_{k+1},q-\tilde m_{k+1}\}\leq q$,
because $1\leq\tilde m_{k+1}\leq q$;
moreover, equality holds only if
$\tilde m_{k+1}=q$.
So $c_k=q=\tilde m_{k+1}$ and $p-q\leq q$.
Since $c_k=q=m_{k+1}$, the index $k$
is of type~B for $m$; as
$m_k=p\neq q=m_{k+1}$, case~(B) of
Theorem~\ref{teo:descripcion}
gives $m_k=m_{k+1}+m_{k+2}$, that is,
$m_{k+2}=p-q=d$.

We keep the notation $\lambda$, $p$, $q$
and $d$ whenever the situation of~(a)
occurs.

(b) Note that $t$ exists, because
$m_{r-1}=\tilde m_{r-1}=c_{r-1}+1$ by
Theorem~\ref{teo:descripcion}.
Let $k>t$ be a position where $m$ and
$\tilde m$ differ. Let $k'\leq k$ be the least index such that
$m$ and $\tilde m$ differ at all positions
$k',k'+1,\ldots,k$.
Then $k'>t$, because $m$ and $\tilde m$
agree at $t$, and they agree at $k'-1$
by the choice of $k'$.
Applying~(a) at $k'$ (after exchanging
$m$ and $\tilde m$ if necessary, so that
$m_{k'}>\tilde m_{k'}$), we get that $m$ and
$\tilde m$ agree at $k'+1$; since they differ
at all positions $k',\ldots,k$, this forces
$k'=k$. Hence they agree at $k-1$ and at $k+1$,
and the last assertion follows from~(a).
\end{proof}

\begin{lemma}\label{lem:tramo-inicial}
Let $m=(m_0,\ldots,m_{r-1},1)$ and
$\tilde m=(\tilde m_0,\ldots,\tilde m_{r-1},1)$
be the multiplicity sequences of two plane
algebroid branches with the same
concentration sequence
$\mathcal{S}=(c_0,\ldots,c_{r-1})$, assume
that $m_0\neq\tilde m_0$, and let $t\geq 1$
be the least index with $m_t=\tilde m_t$.
After possibly exchanging $m$ and $\tilde m$,
exactly one of the following three
situations occurs:
\begin{enumerate}
\item[(I)] $t=1$, the index $0$ is of type~B
for both $m$ and $\tilde m$, and $m_1=\tilde m_1=c_0$;
\item[(II)] $t=2$,
$m=(c_0,c_0,\tilde m_1,c_0-\tilde m_1,\ldots)$
and $\tilde m=(c_0+\tilde m_1,\tilde m_1,\tilde m_1,\ldots)$,
with $\tilde m_1<c_0\leq 2\tilde m_1$;
\item[(III)] $t=3$, $m=(m_0,c_0,m_2,m_2,\ldots)$
and $\tilde m=(c_0+\tilde m_1,\tilde m_1,\tilde m_1,m_2,\ldots)$,
with $m_2=c_0-\tilde m_1$,
$m_2<\tilde m_1\leq 2m_2$
and $m_0\in\{c_0,c_0+m_2\}$.
\end{enumerate}
\end{lemma}

\begin{proof}
Note that $r\geq 2$, since for $r=1$
Theorem~\ref{teo:descripcion}
gives $m=\tilde m=(c_0+1,1)$.
We distinguish cases according to
the types of the index $0$ for $m$
and for $\tilde m$.
Suppose that $0$ is of type~A for both
$m$ and $\tilde m$. Then $c_0>m_1$, $c_0>\tilde m_1$,
$m_0=c_0+m_1$ and $\tilde m_0=c_0+\tilde m_1$.
Lemma~\ref{lem:tipoB}(b), applied with $i=0$
(recall that $r\geq 2$), gives $m_2=m_1$ and
$\tilde m_2=\tilde m_1$, so
Proposition~\ref{proposicion concentracion curva}
at the index $1$ yields
$c_1=\max\{m_2,m_1-m_2\}=m_1$ and,
similarly, $c_1=\tilde m_1$.
Hence $m_1=\tilde m_1$, and therefore
$m_0=c_0+m_1=c_0+\tilde m_1=\tilde m_0$,
contradicting $m_0\neq\tilde m_0$.
If $0$ is of type~B for both
$m$ and $\tilde m$, then, by case~(B) of
Theorem~\ref{teo:descripcion},
$m_1=c_0=\tilde m_1$, that is,
$m_1=\tilde m_1=c_0$; since $m_0\neq\tilde m_0$,
the least index where $m$ and $\tilde m$
agree is $t=1$, and we are in~(I).
Otherwise, $0$ is of type~B for one
of the sequences and of type~A for the
other; exchanging $m$ and $\tilde m$ if
necessary, we may assume without loss of
generality that $0$ is of type~B for $m$
and of type~A for $\tilde m$.
Then $m_1=c_0$, because $0$ is of
type~B for $m$. Since $0$ is of type~A
for $\tilde m$, case~(A) of
Theorem~\ref{teo:descripcion}
gives $\tilde m_1<c_0$ and
$\tilde m_0=c_0+\tilde m_1$. Moreover,
Lemma~\ref{lem:tipoB}(b) applied to
$\tilde m$ gives $\tilde m_2=\tilde m_1$.
Since $\tilde m_1\geq 2$, we have
$\tilde m_2\geq 2$, so $r\geq 3$.
By Proposition~\ref{proposicion concentracion curva},
$c_1=\max\{\tilde m_2,\tilde m_1-\tilde m_2\}
=\tilde m_2=\tilde m_1$
in $\tilde m$ and
$c_1=\max\{m_2,m_1-m_2\}=\max\{m_2,c_0-m_2\}$
in $m$, so $\max\{m_2,c_0-m_2\}=\tilde m_1$.
Hence either $m_2=\tilde m_1$, or
$c_0-m_2=\tilde m_1>m_2$; we consider
these two possibilities separately.

If $m_2=\tilde m_1$, then
$\max\{m_2,c_0-m_2\}=\tilde m_1$ gives
$c_0-\tilde m_1\leq\tilde m_1$, that is,
$c_0\leq 2\tilde m_1$.
Since $m_1=c_0\neq\tilde m_1$
and $m_2=\tilde m_1=\tilde m_2$, we get $t=2$.
As $0$ is of type~B for $m$, case~(B) of
Theorem~\ref{teo:descripcion}
gives $m_0\in\{m_1,m_1+m_2\}=\{c_0,c_0+\tilde m_1\}$,
and $m_0\neq\tilde m_0=c_0+\tilde m_1$
gives $m_0=c_0$.
Moreover, $c_1=\tilde m_1=m_2$, so $1$ is of
type~B for $m$, and case~(B) of
Theorem~\ref{teo:descripcion}
gives $m_1\in\{m_2,m_2+m_3\}$; as
$m_1=c_0>\tilde m_1=m_2$, we get
$m_1=m_2+m_3$, that is,
$m_3=c_0-\tilde m_1$: this is~(II).

If $c_0-m_2=\tilde m_1>m_2$, then
$m_2=c_0-\tilde m_1<\tilde m_1$ and
$c_1=\tilde m_1>m_2$, so $1$ is of
type~A for $m$, and Lemma~\ref{lem:tipoB}(b),
applied with $i=1$ (recall that $r\geq 3$),
gives $m_3=m_2$.
By Proposition~\ref{proposicion concentracion curva}
at the index~$2$, in $m$ we get
$c_2=\max\{m_3,m_2-m_3\}=\max\{m_2,0\}=m_2$,
and in $\tilde m$, using $\tilde m_2=\tilde m_1$,
we get
$c_2=\max\{\tilde m_3,\tilde m_1-\tilde m_3\}$.
Hence $\max\{\tilde m_3,\tilde m_1-\tilde m_3\}=m_2$,
so $\tilde m_3\leq m_2<\tilde m_1$ and
$\tilde m_1-\tilde m_3\leq m_2$; the latter gives
$\tilde m_1\leq\tilde m_3+m_2\leq 2m_2$.
Now $\tilde m_0=c_0+\tilde m_1$, because $0$
is of type~A for $\tilde m$, and
$c_0=\tilde m_1+m_2$, so
$\tilde m_0=2\tilde m_1+m_2$.
Let $(\lambda_0,\ell_0),\ldots,(\lambda_s,\ell_s)$
be the blocks of $\tilde m$.
Since $\tilde m_0>\tilde m_1$, the first block
of $\tilde m$ is
$(\lambda_0,\ell_0)=(2\tilde m_1+m_2,1)$.
Since $\tilde m_1=\tilde m_2>m_2\geq\tilde m_3$,
the second block is
$(\lambda_1,\ell_1)=(\tilde m_1,2)$,
and the third block starts at position~$3$,
so $\lambda_2=\tilde m_3$.
Moreover $1\leq s-1$, because
$\lambda_1=\tilde m_1\geq 2>1=\lambda_s$.
As $0<m_2<\tilde m_1$, $\tilde m_1$ does not
divide $2\tilde m_1+m_2$, so condition~(1) of
Proposition~\ref{prop:mult-to-semigrup}
gives $\tilde m_3=\lambda_2
=\lambda_0-\ell_1\lambda_1
=2\tilde m_1+m_2-2\tilde m_1=m_2$.
Hence $m_1=c_0\neq\tilde m_1$,
$m_2\neq\tilde m_1=\tilde m_2$ and
$m_3=m_2=\tilde m_3$, so $t=3$.
Finally, as $0$ is of type~B for $m$,
case~(B) of
Theorem~\ref{teo:descripcion}
gives $m_0\in\{m_1,m_1+m_2\}=\{c_0,c_0+m_2\}$:
this is~(III).
\end{proof}

\begin{lemma}\label{lem:divisibilidad}
Let $m=(m_0,\ldots,m_{r-1},1)$ and
$\tilde m=(\tilde m_0,\ldots,\tilde m_{r-1},1)$
be the multiplicity sequences of two plane
algebroid branches with the same
concentration sequence
$\mathcal{S}=(c_0,\ldots,c_{r-1})$.
Let $k\geq 1$ be such that
$m_{k-1}=\tilde m_{k-1}=:\lambda$ and
$m_k\neq\tilde m_k$, and, if $k=2$, assume
moreover that $m_0=\tilde m_0$. Then:
\begin{enumerate}
\item[(a)] if $k\geq 3$, then
$m_{k-2}=\tilde m_{k-2}=\lambda$ and
$c_{k-3}=c_{k-2}=\lambda$;
\item[(b)] $\lambda$ divides $m_i$ and
$\tilde m_i$ for every $i<k$.
\end{enumerate}
In particular, $\lambda$ is one of the
integers $e_0,\ldots,e_g$
of~\eqref{notacion:ek-nk}, both for any
branch with multiplicity sequence $m$ and for
any branch with multiplicity sequence
$\tilde m$.
\end{lemma}

\begin{proof}
If $k=1$, then $m_0=\tilde m_0=\lambda$ and
there is nothing to prove, so let $k\geq 2$.
Exchanging $m$ and $\tilde m$ if necessary,
we may assume that $p:=m_k>q:=\tilde m_k$.
By Lemma~\ref{lem:dif-aislada}(a),
$\lambda=p+q$. Since $q=\tilde m_k\geq 1$
is a multiplicity and $q<p$, we get
$0<q<p<p+q=\lambda$.

We first look at the position $k-2$.
Since $m_k=p<\lambda=m_{k-1}$, the index
$k-2$ is not of type~A for $m$: otherwise
Lemma~\ref{lem:tipoB}(b), applied with
$i=k-2$, would give $m_k=m_{k-1}$.
Hence $c_{k-2}=m_{k-1}=\lambda$, and
Lemma~\ref{lem:tipoB}(a) gives
$m_{k-2}\in\{\lambda,\lambda+p\}$.
In the same way,
$\tilde m_{k-2}\in\{\lambda,\lambda+q\}$.
If $k=2$, then $m_0=\tilde m_0$ by
hypothesis, and since
$\lambda+p\neq\lambda+q$, we get
$m_0=\tilde m_0=\lambda$; this proves~(b)
for $k=2$. From now on we assume $k\geq 3$.

(a) Recall that
$\lambda=m_{k-1}=\tilde m_{k-1}=p+q$, with
$p=m_k>q=\tilde m_k$, and that
$m_{k-2}\in\{\lambda,\lambda+p\}$ and
$\tilde m_{k-2}\in\{\lambda,\lambda+q\}$.
We compute $c_{k-3}$ by
Proposition~\ref{proposicion concentracion curva},
first using $m$ and then using $\tilde m$.
Since we do not know $m_{k-3}$, we
distinguish two cases.
\begin{itemize}
\item If $m_{k-2}=\lambda+p$, then
$m_{k-1}=\lambda<m_{k-2}$. If the index
$k-3$ were of type~A for $m$,
Lemma~\ref{lem:tipoB}(b), applied with
$i=k-3$, would give $m_{k-1}=m_{k-2}$, which
is false. Hence $k-3$ is of type~B for $m$,
that is, $c_{k-3}=m_{k-2}=\lambda+p$.
\item If $m_{k-2}=\lambda$, then the sequence
$(m_{k-2},\ldots,m_{r-1},1)$, which is the
multiplicity sequence of the $(k-2)$-th
strict transform of any plane algebroid
branch with multiplicity sequence $m$, begins
with exactly two copies of $\lambda$,
followed by $m_k=p$. By
Lemma~\ref{lem:delante}, applied to this
sequence, with $\ell=2$ and $\rho=p$, and to
$m'=m_{k-3}\geq m_{k-2}$, we get
$m_{k-3}\in\{\lambda,2\lambda,2\lambda+p\}$.
Then
Proposition~\ref{proposicion concentracion curva}
gives $c_{k-3}=\max\{\lambda,m_{k-3}-\lambda\}$,
which is $\lambda$, $\lambda$ or $\lambda+p$,
respectively.
\end{itemize}
Thus, using $m$, we get
$c_{k-3}\in\{\lambda,\lambda+p\}$, and
$c_{k-3}=\lambda+p$ if $m_{k-2}=\lambda+p$.
In the same way, using $\tilde m$, we get
$c_{k-3}\in\{\lambda,\lambda+q\}$, and
$c_{k-3}=\lambda+q$ if
$\tilde m_{k-2}=\lambda+q$.
Since both sequences have the same
$c_{k-3}$ and $p\neq q$, we get
$c_{k-3}=\lambda$. In particular
$m_{k-2}\neq\lambda+p$ and
$\tilde m_{k-2}\neq\lambda+q$, and therefore
$m_{k-2}=\tilde m_{k-2}=\lambda$.
Together with $c_{k-2}=\lambda$, proved
above, this proves part~(a) of the
statement.

(b) We show by descending induction on $i$
that $\lambda$ divides $m_i$ and
$\tilde m_i$ for every $i\leq k-1$.
For $i=k-1$ and $i=k-2$ this holds, since
$m_i=\tilde m_i=\lambda$ by the definition of
$\lambda$ and by~(a).
Let $i\leq k-3$, and assume that $\lambda$
divides $m_j$ and $\tilde m_j$ for
$i<j\leq k-1$. We first work with $m$.
Write $\ell=\ell_{i+1}(m)$ and
$\rho=\rho_{i+1}(m)$.
The sequence $(m_{i+1},\ldots,m_{r-1},1)$ is
the multiplicity sequence of the
$(i+1)$-th strict transform of any plane
algebroid branch with multiplicity sequence
$m$, and $(m_i,\ldots,m_{r-1},1)$ is that of
its $i$-th strict transform; hence
Lemma~\ref{lem:delante}, applied to the
first sequence and to
$m'=m_i\geq m_{i+1}$, gives
\[
m_i\in\{m_{i+1},2m_{i+1},\ldots,
\ell m_{i+1},\ell m_{i+1}+\rho\}.
\]
Since $m_{i+1}\geq m_{k-1}=\lambda>p=m_k$,
the $\ell$ consecutive terms equal to
$m_{i+1}$ starting with $m_{i+1}$ end at a
position $i+\ell\leq k-1$. We distinguish
two cases.
\begin{itemize}
\item If $i+\ell<k-1$, then
$\rho=m_{i+\ell+1}$ with
$i+\ell+1\leq k-1$, so $\lambda$ divides
$m_{i+1}$ and $\rho$ by the induction
hypothesis. Hence every element of the set
above is a multiple of $\lambda$, so
$\lambda$ divides $m_i$, and also
$c_i=\max\{m_{i+1},m_i-m_{i+1}\}$.
\item If $i+\ell=k-1$, then
$m_{i+1}=\cdots=m_{k-1}=\lambda$; moreover
$\ell\geq 2$, because $i+1\leq k-2$, and
$\rho=m_k=p$. So the set above is
$\{\lambda,2\lambda,\ldots,\ell\lambda,
\ell\lambda+p\}$.
If $m_i=j\lambda$ for some
$1\leq j\leq\ell$, then $\lambda$ divides
$m_i$ and
$c_i=\max\{\lambda,(j-1)\lambda\}$.
If $m_i=\ell\lambda+p$, then
Proposition~\ref{proposicion concentracion curva}
gives
$c_i=\max\{\lambda,(\ell-1)\lambda+p\}
=(\ell-1)\lambda+p$.
\end{itemize}
Hence, using $m$, we get
$c_i\equiv 0\pmod\lambda$ if $\lambda$
divides $m_i$, and $c_i\equiv p\pmod\lambda$
otherwise.
In the same way, using $\tilde m$, we get
$c_i\equiv 0\pmod\lambda$ if $\lambda$
divides $\tilde m_i$, and
$c_i\equiv q\pmod\lambda$ otherwise.
Since $0<q<p<\lambda$, the residues $0$, $q$
and $p$ modulo $\lambda$ are pairwise
distinct. As both sequences have the same
$c_i$, its residue modulo $\lambda$ must be
the same in the two computations, which is
only possible if it is $0$ in both. Hence
$\lambda$ divides both $m_i$ and
$\tilde m_i$. This proves part~(b) of the
statement.

Finally, we prove the last assertion of
the statement. In $m$ the block of value
$\lambda$ ends at the position $k-1$,
because $m_k<\lambda$.
If it is the first block, then
$\lambda=m_0=e_0$.
Otherwise, it is preceded by a block whose
value is some $m_i$ with $i<k-1$, and
$\lambda$ divides $m_i$ by~(b); by the
comment after
Lemma~\ref{lem:bloques-divide},
$\lambda$ is one of $e_1,\ldots,e_g$.
The same argument applies to $\tilde m$.
\end{proof}

We now show that $\mathcal{M}(\mathcal{S})$
is, in a natural way, a Cartesian product:
an element of $\mathcal{M}(\mathcal{S})$ is
determined by its first three terms and by
its terms at certain positions, which we
call free positions, and these terms can be
chosen independently, each of them among two
values (Theorem~\ref{teo:producto}).
This allows us to compute explicitly the
number $N(\mathcal{S})$ of elements of
$\mathcal{M}(\mathcal{S})$ in terms of
$\mathcal{S}$. We first note that each term
of a sequence in $\mathcal{M}(\mathcal{S})$
can only take two values.

\begin{lemma}\label{lem:dos-valores}
Let $m=(m_0,\ldots,m_{r-1},1)$ be the
multiplicity sequence of a plane algebroid
branch with concentration sequence
$\mathcal{S}=(c_0,\ldots,c_{r-1})$.
Then $m_k\in\{c_{k-1},c_k\}$ for every
$1\leq k\leq r-1$.
\end{lemma}

\begin{proof}
By Proposition~\ref{proposicion concentracion curva},
$c_{k-1}=\max\{m_k,m_{k-1}-m_k\}\geq m_k$.
If $c_{k-1}=m_k$, we are done.
Otherwise $c_{k-1}>m_k$, and
Lemma~\ref{lem:tipoB}(b), applied with
$i=k-1$, gives $m_{k+1}=m_k$; then
Proposition~\ref{proposicion concentracion curva}
gives
$c_k=\max\{m_{k+1},m_k-m_{k+1}\}=m_k$.
\end{proof}

\begin{definition}\label{def:libre}
Let $\mathcal{S}=(c_0,\ldots,c_{r-1})$ be a
sequence of positive integers.
A position $k$ with $3\leq k\leq r-2$ is
\emph{free} for $\mathcal{S}$ if, writing
$\lambda=c_{k-1}+c_k$, the following three
conditions hold:
\begin{enumerate}
\item[(F1)] $c_k<c_{k-1}\leq 2c_k$;
\item[(F2)] $c_{k-3}=c_{k-2}=\lambda$;
\item[(F3)] $\lambda$ divides $c_j$ for every
$0\leq j\leq k-2$.
\end{enumerate}
\end{definition}

We denote by $\mathfrak{f}(\mathcal{S})$ the
number of free positions of $\mathcal{S}$,
and by $\mathcal{T}(\mathcal{S})$ the set of
triples $(m_0,m_1,m_2)$ formed by the first
three terms of the elements $m$ of
$\mathcal{M}(\mathcal{S})$. Note that
$\mathfrak{f}(\mathcal{S})$ can be computed
directly from $\mathcal{S}$.

\begin{theorem}\label{teo:producto}
Let $\mathcal{S}=(c_0,\ldots,c_{r-1})$ be a
sequence of positive integers with $r\geq 3$
and $\mathcal{M}(\mathcal{S})\neq\emptyset$.
\begin{enumerate}
\item[(a)] If $3\leq k\leq r-1$ is not a free
position of $\mathcal{S}$, then all the
elements of $\mathcal{M}(\mathcal{S})$ have
the same term at the position $k$.
\item[(b)] Let $k$ be a free position of
$\mathcal{S}$, and write $p=c_{k-1}$,
$q=c_k$ and $\lambda=p+q$.
Every $m\in\mathcal{M}(\mathcal{S})$ satisfies
\[
m_{k-2}=m_{k-1}=\lambda,\quad
m_k\in\{p,q\},\quad
m_{k+1}=q,\quad
m_{k+2}=p-q,
\]
and the sequence obtained from $m$ by
replacing $m_k$ with the other element of
$\{p,q\}$ belongs to $\mathcal{M}(\mathcal{S})$.
\item[(c)] The set $\mathcal{T}(\mathcal{S})$
is contained in the set formed by the five
triples
\[
\begin{aligned}
&(c_0,c_0,c_1),\quad (c_0+c_1,c_0,c_1),\quad
(c_0,c_0,c_2),\\
&(c_0+c_2,c_0,c_2),\quad (c_0+c_1,c_1,c_1),
\end{aligned}
\]
and $\#\mathcal{T}(\mathcal{S})$ is $1$, $2$,
$3$ or $5$.
\item[(d)] The map that sends
$m\in\mathcal{M}(\mathcal{S})$ to
$(m_0,m_1,m_2)$ and to its terms $m_k$ at the
free positions $k$ of $\mathcal{S}$ is a
bijection from $\mathcal{M}(\mathcal{S})$
onto
\[
\mathcal{T}(\mathcal{S})\times
\prod_{k\ \mathrm{free}}\{c_{k-1},c_k\},
\]
the Cartesian product of
$\mathcal{T}(\mathcal{S})$ and of the sets
$\{c_{k-1},c_k\}$, where $k$ runs over the
free positions of $\mathcal{S}$ (if there
are no free positions, this is just
$\mathcal{T}(\mathcal{S})$). In particular,
\[
N(\mathcal{S})=\#\mathcal{T}(\mathcal{S})\,
2^{\mathfrak{f}(\mathcal{S})}.
\]
\end{enumerate}
\end{theorem}
\begin{proof}
(a) Let $m,\tilde m\in\mathcal{M}(\mathcal{S})$
and $3\leq k\leq r-1$ with
$m_k\neq\tilde m_k$; we show that $k$ is free.
Since $m$ and $\tilde m$ belong to
$\mathcal{M}(\mathcal{S})$, both satisfy the
conditions of
Theorem~\ref{teo:descripcion}; the first
one gives $m_{r-1}=c_{r-1}+1=\tilde m_{r-1}$.
As $m$ and $\tilde m$ differ at $k$, we get
$k\neq r-1$, that is, $k\leq r-2$.
Let $t$ be the least index with
$m_t=\tilde m_t$; it exists
because $m_{r-1}=\tilde m_{r-1}$.
If $m_0=\tilde m_0$, then $t=0$.
Otherwise, $m$ and $\tilde m$
satisfy the hypotheses of
Lemma~\ref{lem:tramo-inicial}: they have
the same concentration sequence $\mathcal{S}$
and $m_0\neq\tilde m_0$. In the situations
(I), (II) and (III) of this lemma we have
$t=1$, $t=2$ and $t=3$, respectively; hence
$t\leq 3$.
Since $m$ and $\tilde m$ agree at
$t$ and differ at $k$, we have $k\neq t$, and
$k\geq 3\geq t$; hence $k>t$.
Thus $k$ is a position after $t$
where $m$ and $\tilde m$ differ, and
Lemma~\ref{lem:dif-aislada}(b) gives that
$m$ and $\tilde m$ agree at $k-1$.
Exchanging $m$ and $\tilde m$ if necessary,
Lemma~\ref{lem:dif-aislada}(a) gives
$p:=m_k>q:=\tilde m_k$,
$\lambda:=m_{k-1}=p+q$,
$\tilde m_{k+1}=q$ and $p-q\leq q$.
By
Proposition~\ref{proposicion concentracion curva},
applied to the strict transforms of branches
with multiplicity sequences $m$ and
$\tilde m$, both with concentration sequence
$\mathcal{S}$,
\[
c_i=\max\{m_{i+1},m_i-m_{i+1}\}
=\max\{\tilde m_{i+1},\tilde m_i-\tilde m_{i+1}\}
\quad\text{for } 0\leq i\leq r-1.
\]
For $i=k-1$ in $m$, where $m_{k-1}=\lambda$
and $m_k=p$, this formula reads
$c_{k-1}=\max\{p,\lambda-p\}=\max\{p,q\}=p$,
and for $i=k$ in $\tilde m$, where
$\tilde m_k=\tilde m_{k+1}=q$, it yields
$c_k=\max\{q,0\}=q$.
Since $p-q\leq q$, we have
$c_k=q<p=c_{k-1}\leq 2q=2c_k$, which
is condition~(F1) of
Definition~\ref{def:libre}, and
$\lambda=p+q=c_{k-1}+c_k$.
Now $m$ and $\tilde m$ satisfy the
hypotheses of
Lemma~\ref{lem:divisibilidad}:
$m_{k-1}=\tilde m_{k-1}=\lambda$,
$m_k\neq\tilde m_k$, and $k\geq 3$ (so the
extra hypothesis for $k=2$ is not needed).
Part~(a) of this lemma gives
$c_{k-3}=c_{k-2}=\lambda$, which is~(F2).
Part~(b) gives that $\lambda$ divides $m_i$
for every $i\leq k-1$. Let $0\leq j\leq k-2$.
Then $j+1\leq k-1$, so $\lambda$ divides
$m_j$ and $m_{j+1}$, hence also
$m_j-m_{j+1}$, and therefore it divides
$c_j=\max\{m_{j+1},m_j-m_{j+1}\}$, which is
one of these two numbers. This is~(F3).
Finally, $3\leq k\leq r-2$, as shown at the
beginning, so $k$ is a free position of
$\mathcal{S}$.

(b) Let $m\in\mathcal{M}(\mathcal{S})$.
Since $k$ is a free position,
condition~(F2) of
Definition~\ref{def:libre} gives
$c_{k-3}=c_{k-2}=\lambda$, and, by
definition, $c_{k-1}=p$ and $c_k=q$.
Since $1\leq k-2<k\leq r-2$,
Lemma~\ref{lem:dos-valores} applies at
the positions $k-2$, $k-1$ and $k$, and gives
$m_{k-2}\in\{c_{k-3},c_{k-2}\}=\{\lambda\}$,
so $m_{k-2}=\lambda$; moreover,
$m_{k-1}\in\{c_{k-2},c_{k-1}\}=\{\lambda,p\}$
and $m_k\in\{c_{k-1},c_k\}=\{p,q\}$.
If $m_{k-1}=p$, then
Proposition~\ref{proposicion concentracion curva}
gives
$c_{k-2}=\max\{m_{k-1},m_{k-2}-m_{k-1}\}
=\max\{p,\lambda-p\}=\max\{p,q\}=p$,
which contradicts $c_{k-2}=\lambda>p$.
Hence $m_{k-1}=\lambda$.

If $m_k=p$, we use
Lemma~\ref{lem:tipoB}(a) twice, in the
following form: if $c_i=m_{i+1}$ and
$m_i\neq m_{i+1}$, then
$m_i=m_{i+1}+m_{i+2}$. For $i=k-1$ we have
$c_{k-1}=p=m_k$ and $m_{k-1}=\lambda\neq p$,
so $m_{k+1}=\lambda-p=q$. For $i=k$ we now
have $c_k=q=m_{k+1}$ and $m_k=p\neq q$, so
$m_{k+2}=p-q$.

If $m_k=q$, then $c_{k-1}=p>q=m_k$,
so case~(A) of
Theorem~\ref{teo:descripcion} holds at
$k-1$, and Lemma~\ref{lem:tipoB}(b),
applied with $i=k-1$, gives $m_{k+1}=q$.
Hence $\ell=\ell_k(m)\geq 2$, and
$\rho=\rho_k(m)$ satisfies $0<\rho<q$. By~(A),
either $q$ divides $p$ and $p\leq(\ell-1)q$,
or $p=(\ell-1)q+\rho$. Since $q<p\leq 2q$
by~(F1), in the first case $p=2q$ and
$\ell\geq 3$, so $m_{k+2}=q=p-q$; in the
second case $(\ell-1)q<p\leq 2q$ forces
$\ell=2$, so $m_{k+2}=\rho=p-q$.

Now let $m'$ be the sequence obtained from
$m$ by replacing $m_k$ with the other element
of $\{p,q\}$. We check that $m'$ satisfies the
conditions of
Theorem~\ref{teo:descripcion}.
Since $k\leq r-2$, we have
$m'_{r-1}=m_{r-1}=c_{r-1}+1$.
Let $0\leq i\leq r-2$.
The condition at $i$ only involves $m'_i$,
$m'_{i+1}$, $m'_{i+2}$, $\ell_{i+1}(m')$ and
$\rho_{i+1}(m')$, so for $i\geq k+1$ it is the
same for $m'$ as for $m$.
Recall that $m'_{k-2}=m'_{k-1}=\lambda$,
$m'_{k+1}=q$ and $m'_{k+2}=p-q$.
\begin{itemize}
\item At $i=k$: $c_k=q=m'_{k+1}$ and
$m'_k\in\{p,q\}=\{q,q+m'_{k+2}\}$,
so~(B) holds.
\item At $i=k-1$: if $m'_k=p$, then
$c_{k-1}=m'_k$ and
$m'_{k-1}=\lambda=m'_k+m'_{k+1}$, so~(B)
holds. If $m'_k=q$, then $c_{k-1}=p>m'_k$,
$m'_{k-1}=\lambda=c_{k-1}+m'_k$, and, as
$m'_k=m'_{k+1}=q$ and $m'_{k+2}=p-q$, either
$p<2q$, $\ell_k(m')=2$ and $p=q+\rho_k(m')$,
or $p=2q$ and $\ell_k(m')\geq 3$, so that
$q$ divides $p$ and $p\leq(\ell_k(m')-1)q$;
in both cases~(A) holds.
\item At $i=k-2$:
$c_{k-2}=\lambda=m'_{k-1}$ and
$m'_{k-2}=\lambda$, so~(B) holds.
\item At $i\leq k-3$: the terms $m'_i$,
$m'_{i+1}$, $m'_{i+2}$ are those of $m$, so
(B) at $i$ passes from $m$ to $m'$. Assume
that (A) holds at $i$ for $m$. Changing the
term at $k$ can only change $\rho_{i+1}$
(not $\ell_{i+1}$, as $m_k,m'_k<\lambda$),
and only if $m_{i+1}=\cdots=m_{k-1}=\lambda$.
In that case $\lambda$ divides $c_i$
by~(F3), so
$c_i\neq(\ell_{i+1}(m)-1)\lambda+\rho_{i+1}(m)$,
because $0<\rho_{i+1}(m)<\lambda$. Hence (A)
holds for $m$ through its first alternative,
which does not involve $\rho_{i+1}$, and so
it holds for $m'$.
\end{itemize}
By Theorem~\ref{teo:descripcion},
$m'\in\mathcal{M}(\mathcal{S})$.

(c) Let $m\in\mathcal{M}(\mathcal{S})$.
By Lemma~\ref{lem:dos-valores},
$m_1\in\{c_0,c_1\}$ and $m_2\in\{c_1,c_2\}$.
If $m_1\neq c_0$, then $m_1=c_1$ and
$c_0>m_1$, because $c_0\geq m_1$ by
Proposition~\ref{proposicion concentracion curva};
so~(A) holds at $0$, $m_0=c_0+c_1$, and
Lemma~\ref{lem:tipoB}(b) gives
$m_2=m_1=c_1$. Thus
$(m_0,m_1,m_2)=(c_0+c_1,c_1,c_1)$.
If $m_1=c_0$, then~(B) holds at $0$, and
$m_0\in\{c_0,c_0+m_2\}$ with
$m_2\in\{c_1,c_2\}$. This proves the first
assertion.

Moreover, if (B) holds at $0$ for $m$, then
replacing $m_0$ with the other element of
$\{m_1,m_1+m_2\}$ gives again an element of
$\mathcal{M}(\mathcal{S})$, by
Theorem~\ref{teo:descripcion}, because
$m_0$ appears in no other condition.
Hence
$\#\mathcal{T}(\mathcal{S})=2b+a$, where
$b$ is the number of $x\in\{c_1,c_2\}$
with $(c_0,c_0,x)\in\mathcal{T}(\mathcal{S})$,
and $a\in\{0,1\}$ is the number of
elements $(m_0,m_1,m_2)$ of
$\mathcal{T}(\mathcal{S})$ with
$m_1\neq c_0$.
Thus $b\in\{0,1,2\}$ and
$a\in\{0,1\}$. If $b=0$, then
$a=1$, because
$\mathcal{T}(\mathcal{S})\neq\emptyset$ (as
$\mathcal{M}(\mathcal{S})\neq\emptyset$), so
$\#\mathcal{T}(\mathcal{S})=1$. If $b=1$,
then $\#\mathcal{T}(\mathcal{S})\in\{2,3\}$.
If $b=2$, then
$\#\mathcal{T}(\mathcal{S})\in\{4,5\}$, and it
remains to exclude the value $4$, that is,
to show that $a=1$.

Assume $b=2$. Then $c_1\neq c_2$, and
there are $m,\tilde m\in\mathcal{M}(\mathcal{S})$
with $m_1=\tilde m_1=c_0$ and
$\{m_2,\tilde m_2\}=\{c_1,c_2\}$; say
$m_2>\tilde m_2$.
By Lemma~\ref{lem:dif-aislada}(a),
applied with $k=2$, we get
$c_0=m_2+\tilde m_2$ and
$m_3=\tilde m_2$, and by
Proposition~\ref{proposicion concentracion curva},
$c_1=\max\{m_2,c_0-m_2\}=m_2$.
Hence $m_2=c_1$, $\tilde m_2=c_2$,
$c_0=c_1+c_2$ and $m_3=c_2\neq c_1$.
Let
\[
m''=(c_0+c_1,c_1,c_1,m_3,\ldots,m_{r-1},1).
\]
Since $m''_j=m_j$ for $j\geq 2$, the
conditions of
Theorem~\ref{teo:descripcion} at the
indices $i\geq 2$ are the same for $m''$ as
for $m$, and, since $r\geq 3$,
$m''_{r-1}=m_{r-1}=c_{r-1}+1$.
At $i=1$, $c_1=m''_2=m''_1$, so~(B) holds.
At $i=0$, $c_0=c_1+c_2>c_1=m''_1$ and
$m''_0=c_0+m''_1$; moreover
$\ell_1(m'')=2$ and $\rho_1(m'')=m_3=c_2$, so
$c_0=(\ell_1(m'')-1)c_1+\rho_1(m'')$ and~(A)
holds.
Hence $m''\in\mathcal{M}(\mathcal{S})$, and
$a=1$, so
$\#\mathcal{T}(\mathcal{S})=5$.

(d) The map is well defined by~(b)
and~(c).
It is injective: if two elements of
$\mathcal{M}(\mathcal{S})$ have the same
image, they agree at the positions $0$, $1$,
$2$ and at the free positions, and by~(a)
also at the other positions $k\geq 3$.
It is surjective: given
$\tau\in\mathcal{T}(\mathcal{S})$ and, for
each free position $k$, an element
$x_k\in\{c_{k-1},c_k\}$, take
$m\in\mathcal{M}(\mathcal{S})$ with
$(m_0,m_1,m_2)=\tau$, and, for each free
position $k$ with $m_k\neq x_k$, replace
$m_k$ by $x_k$. By~(b), each replacement
gives an element of $\mathcal{M}(\mathcal{S})$,
and it does not change any other term.
The formula for $N(\mathcal{S})$ follows.
\end{proof}
\begin{examples}\label{ex:producto}
\leavevmode
\begin{enumerate}[label=(\arabic*),leftmargin=*]
\item Let $\mathcal{S}=(15,10,5,5,3,2,1)$, so that
$r=7$ and the possible free positions are
$k=3,4,5$.
For $k=3$, condition~(F1) of
Definition~\ref{def:libre} fails, since
$c_3=c_2=5$.
For $k=4$, condition~(F1) holds, since
$c_4=3<c_3=5\leq 6$, but (F2) fails, since
$\lambda=c_3+c_4=8$ and $c_2=5\neq 8$.
For $k=5$, we have $\lambda=c_4+c_5=5$, and
the three conditions hold:
$c_5=2<c_4=3\leq 4$, $c_2=c_3=5$, and $5$
divides $c_0,c_1,c_2,c_3$, which are $15$,
$10$, $5$ and $5$.
Hence $5$ is the only free position, and
$\mathfrak{f}(\mathcal{S})=1$.
By Theorem~\ref{teo:producto}(b), with
$p=c_4=3$ and $q=c_5=2$, every
$m\in\mathcal{M}(\mathcal{S})$ satisfies
$m_3=m_4=5$, $m_5\in\{3,2\}$, $m_6=2$ and
$m_7=1$.
The five triples of
Theorem~\ref{teo:producto}(c) are
$(15,15,10)$, $(25,15,10)$, $(15,15,5)$,
$(20,15,5)$ and $(25,10,10)$, and the
algorithm of Corollary~\ref{cor:algoritmo}
shows that all of them belong to
$\mathcal{T}(\mathcal{S})$, so
$\#\mathcal{T}(\mathcal{S})=5$.
Therefore, by
Theorem~\ref{teo:producto}(d),
$\mathcal{M}(\mathcal{S})$ consists of the
$5\cdot 2=10$ sequences
\[
(m_0,m_1,m_2,5,5,x,2,1),
\]
where $(m_0,m_1,m_2)$ runs over
$\mathcal{T}(\mathcal{S})$ and $x\in\{3,2\}$,
chosen independently.

\item Let
$\mathcal{S}=(90,45,45,45,30,15,15,10,5,5,3,2,1)$,
so that $r=13$ and the possible free
positions are $3\leq k\leq 11$.
The position $k=5$ is free, with
$\lambda=c_4+c_5=45$: indeed,
$c_5=15<c_4=30\leq 30$, $c_2=c_3=45$, and
$45$ divides $c_0,\ldots,c_3$, which are
$90$, $45$, $45$ and $45$.
In the same way, $k=8$ is free, with
$\lambda=c_7+c_8=15$, and $k=11$ is free,
with $\lambda=c_{10}+c_{11}=5$, and one
checks that no other position satisfies
(F1), (F2) and (F3). Hence
$\mathfrak{f}(\mathcal{S})=3$.
Since $c_1=c_2=45$, the five triples of
Theorem~\ref{teo:producto}(c) reduce to
$(90,90,45)$, $(135,90,45)$ and
$(135,45,45)$, and all of them belong to
$\mathcal{T}(\mathcal{S})$. Therefore
$N(\mathcal{S})=3\cdot 2^3=24$.
\end{enumerate}
\end{examples}

Theorem~\ref{teo:producto} gives the
exact value of $N(\mathcal{S})$, the number
of multiplicity sequences of plane algebroid
branches with concentration sequence
$\mathcal{S}$, but
computing it requires knowing
$\mathcal{S}$: one has to determine
$\#\mathcal{T}(\mathcal{S})$ and the free positions of
$\mathcal{S}$. The following consequence
gives instead a bound that only depends on
the length $r$ of $\mathcal{S}$, that is,
on the number of blowups needed to reach
the first smooth strict transform. It shows
that the number of non-equisingular
branches sharing a concentration sequence
of length $r$ grows at most like $2^{r/3}$.

\begin{corollary}\label{cor:cota-r}
Let $\mathcal{S}=(c_0,\ldots,c_{r-1})$ be a
sequence of positive integers with
$r\geq 3$. Then
\[
N(\mathcal{S})\leq
5\cdot 2^{\lfloor (r-2)/3\rfloor}.
\]
\end{corollary}

\begin{proof}
If $\mathcal{M}(\mathcal{S})=\emptyset$
there is nothing to prove, so assume
$\mathcal{M}(\mathcal{S})\neq\emptyset$.
By Theorem~\ref{teo:producto},
$N(\mathcal{S})=\#\mathcal{T}(\mathcal{S})\,
2^{\mathfrak{f}(\mathcal{S})}$ with
$\#\mathcal{T}(\mathcal{S})\leq 5$, so it suffices to
show that $\mathfrak{f}(\mathcal{S})\leq
\lfloor (r-2)/3\rfloor$.
Let $k<k'$ be two free positions of
$\mathcal{S}$.
If $k'=k+1$, then condition~(F2) of
Definition~\ref{def:libre} for $k'$ gives
$c_{k-2}=c_{k-1}$, while (F2) for $k$ gives
$c_{k-2}=c_{k-1}+c_k>c_{k-1}$, a
contradiction.
If $k'=k+2$, then (F2) for $k'$ gives
$c_{k-1}=c_k$, while (F1) for $k$ gives
$c_k<c_{k-1}$, a contradiction.
Hence $k'\geq k+3$.
Since the free positions lie in
$\{3,\ldots,r-2\}$ and any two of them
differ by at least $3$, there are at most
$\lfloor (r-5)/3\rfloor+1
=\lfloor (r-2)/3\rfloor$ of them.
\end{proof}

The bound of
Corollary~\ref{cor:cota-r} is attained,
for instance, for $r=4$, $7$ and $10$:
for $\mathcal{S}=(7,4,3,2)$ we have
$N(\mathcal{S})=5$
(Example~\ref{ex:arbol}); for
$\mathcal{S}=(15,10,5,5,3,2,1)$ we have
$N(\mathcal{S})=10$
(Examples~\ref{ex:producto}(1)); and for
$\mathcal{S}=(45,30,15,15,10,5,5,3,2,1)$ we
have $\#\mathcal{T}(\mathcal{S})=5$, the free positions
are $5$ and $8$, and $N(\mathcal{S})=20$.

We finish this section by characterizing the
sequences of positive integers that are the
concentration sequence of some plane branch. We
say that $\mathcal{S}$ is a \emph{concentration
sequence} if it is the concentration sequence of
some plane branch, that is, if
$\mathcal{M}(\mathcal{S})\neq\emptyset$. In
principle, this can be decided by computing
$\mathcal{M}(\mathcal{S})$ with the algorithm of
Corollary~\ref{cor:algoritmo}. The following
remark, which is Theorem~\ref{teo:descripcion} read
at the position~$0$, describes how a term can be
added at the beginning of a concentration sequence.

\begin{remark}\label{rem:delante-conc}
Let $\mathcal{S}=(c_0,\ldots,c_{r-1})$ be a
concentration sequence, let
$m=(m_0,\ldots,m_{r-1},1)\in\mathcal{M}(\mathcal{S})$,
and write $\ell=\ell_0(m)$ and $\rho=\rho_0(m)$. We
say that a positive integer $n$ \emph{can be put in
front of $\mathcal{S}$ through~$m$} if some element
of $\mathcal{M}((n,c_0,\ldots,c_{r-1}))$ ends with
$m_0,\ldots,m_{r-1},1$. By
Theorem~\ref{teo:descripcion}, this happens if and
only if
\[
n=m_0,\quad\text{or}\quad \ell\geq 2\ \text{and}\
n\in\{jm_0 : 2\leq j\leq\ell-1\}
\cup\{(\ell-1)m_0+\rho\}.
\]
Indeed, the conditions of
Theorem~\ref{teo:descripcion} for a sequence
$(m',m_0,\ldots,m_{r-1},1)$ at the positions
$i\geq 1$ and on the last term are those of $m$ at
the positions $i-1$, which hold because
$m\in\mathcal{M}(\mathcal{S})$. At the
position~$0$, case~(B) gives $n=m_0$, with
$m'=m_0$ or $m'=m_0+m_1$, and case~(A) gives
$n>m_0$ and $m'=n+m_0$, with either $m_0$ divides
$n$ and $n\leq(\ell-1)m_0$, or
$n=(\ell-1)m_0+\rho$; the latter requires
$\ell\geq 2$, since $0<\rho<m_0<n$. In particular,
taking $m'=m_0$ in case~(B), we get
$(m_0,m_0,\ldots,m_{r-1},1)\in
\mathcal{M}((m_0,c_0,\ldots,c_{r-1}))$.

Since removing the first term of an element of
$\mathcal{M}((n,c_0,\ldots,c_{r-1}))$ gives an
element of $\mathcal{M}(\mathcal{S})$, the sequence
$(n,c_0,\ldots,c_{r-1})$ is a concentration
sequence if and only if $n$ can be put in front of
$\mathcal{S}$ through some
$m\in\mathcal{M}(\mathcal{S})$. Finally, by
Theorem~\ref{teo:descripcion},
$\mathcal{M}((c_0))=\{(c_0+1,1)\}$ for every
$c_0\geq1$, so every sequence $(c_0)$ is a
concentration sequence; as $\ell_0(c_0+1,1)=1$, the
sequence $(n,c_0)$ is a concentration sequence if
and only if $n=c_0+1$.
\end{remark}

To characterize the concentration sequences we
also need the following consequence of
Theorem~\ref{teo:producto}(b).

\begin{lemma}\label{lem:cambio-posicion}
Let $\mathcal{S}=(c_0,\ldots,c_{r-1})$ be a
sequence of positive integers, let
$1\leq k\leq r-2$ be such that
$c_k<c_{k-1}\leq 2c_k$, and let
$m\in\mathcal{M}(\mathcal{S})$ be such that
$m_0=\cdots=m_{k-1}=c_{k-1}+c_k$. Then the sequence
obtained from $m$ by replacing $m_k$ with the other
element of $\{c_{k-1},c_k\}$ belongs to
$\mathcal{M}(\mathcal{S})$.
\end{lemma}

\begin{proof}
Write $\lambda=c_{k-1}+c_k$. For $0\leq i\leq k-2$
we have $m_i=m_{i+1}=\lambda$, so
$c_i=\max\{m_{i+1},m_i-m_{i+1}\}=\lambda$
(Proposition~\ref{proposicion concentracion curva}).
By Remark~\ref{rem:delante-conc}, applied twice
with $n=\lambda=m_0$, the sequence
$(\lambda,\lambda,m_0,\ldots,m_{r-1},1)$ belongs to
$\mathcal{M}((\lambda,\lambda,c_0,\ldots,c_{r-1}))$.
In $(\lambda,\lambda,c_0,\ldots,c_{r-1})$, the
terms at the positions $0,\ldots,k$ are equal to
$\lambda$, those at the positions $k+1$ and $k+2$
are $c_{k-1}$ and $c_k$, and $3\leq k+2\leq r$;
hence the position $k+2$ is free
(Definition~\ref{def:libre}). It remains to apply
Theorem~\ref{teo:producto}(b) at this position,
where the term is $m_k$, and to remove the first
two terms.
\end{proof}

The next proposition characterizes the concentration sequences from
right to left, by describing the integers that can be put in front of
a concentration sequence.

\begin{proposition}\label{prop:imagen}
Let $\mathcal{S}=(c_0,\ldots,c_{r-1})$ be a
concentration sequence with $r\geq2$, put $c_r=1$,
and let $\ell$ be the number of initial terms of
$\mathcal{S}$ equal to $c_0$. Then
$(n,c_0,\ldots,c_{r-1})$ is a concentration
sequence if and only if $n$ is one of the following
integers:
\begin{enumerate}
\item[(i)] $c_0,2c_0,\ldots,\ell c_0$, if
$c_0\leq 2c_\ell$;
\item[(ii)] $\ell c_0+c_\ell$, if
$c_\ell\leq 2c_{\ell+1}$ and $\ell\leq r-2$;
\item[(iii)] $\ell c_0+c_{\ell+1}$, if
$c_0=c_\ell+c_{\ell+1}$ and $\ell\neq r-2$.
\end{enumerate}
\end{proposition}

\begin{proof}
Let $m\in\mathcal{M}(\mathcal{S})$. Recall that
$c_i=\max\{m_{i+1},m_i-m_{i+1}\}$ for
$0\leq i\leq r-1$, where $m_r=1$
(Proposition~\ref{proposicion concentracion curva}).
If $0$ is of type~A for $m$, then $c_0>m_1$ and
$m_0=c_0+m_1$, and Lemmas~\ref{lem:tipoB}(b)
and~\ref{lem:dos-valores} give $m_2=m_1=c_1<c_0$.
Hence $\ell=1$; moreover $r\geq3$, as
$m_2=m_1\geq2$, and $1$ is of type~B for $m$, so
$c_1\leq 2c_2$ by Remark~\ref{rem:algoritmo}(b). As
$m_0>m_1$, the only integer that can be put in
front of $\mathcal{S}$ through $m$ is
$m_0=c_0+c_1=\ell c_0+c_\ell$, which is given
by~(ii).

If $0$ is of type~B for $m$, then $m_1=c_0$, and
replacing $m_0$ with the other element of
$\{m_1,m_1+m_2\}$ gives again an element of
$\mathcal{M}(\mathcal{S})$, by
Theorem~\ref{teo:descripcion}. If $m_i=c_0$ for
some $1\leq i\leq\ell-1$, then
$c_0=c_i=\max\{m_{i+1},c_0-m_{i+1}\}$ forces
$m_{i+1}=c_0$. Hence $m_1=\cdots=m_\ell=c_0$, so
$\ell\leq r-1$, as $m_r=1$, and
$c_\ell=\max\{m_{\ell+1},c_0-m_{\ell+1}\}<c_0$
gives $m_{\ell+1}<c_0$. Moreover, $\ell-1$ is of
type~B for $m$, so $c_0\leq 2c_\ell$ by
Remark~\ref{rem:algoritmo}(b). If $m_0=c_0$, then
$\ell_0(m)=\ell+1$ and $\rho_0(m)=m_{\ell+1}$, so
the integers that can be put in front through $m$
are $c_0,2c_0,\ldots,\ell c_0$ and
$\ell c_0+m_{\ell+1}$. If $m_0=c_0+m_2$, the only
one is $c_0+m_2$, which is $2c_0$ if $\ell\geq2$
and $\ell c_0+m_{\ell+1}$ if $\ell=1$, so it can
also be put in front through the element obtained
by replacing $m_0$ with $c_0$.

Assume now that $m_0=m_1=c_0$. If $\ell=r-1$, then
$m_{\ell+1}=m_r=1=c_{\ell+1}$ and
$c_0=m_{r-1}=c_{r-1}+1=c_\ell+c_{\ell+1}$
(Theorem~\ref{teo:descripcion}), so
$\ell c_0+m_{\ell+1}$ is given by~(iii). If
$\ell\leq r-2$, then
$m_{\ell+1}\in\{c_\ell,c_{\ell+1}\}$ by
Lemma~\ref{lem:dos-valores}. If
$m_{\ell+1}=c_\ell$, then $\ell$ is of type~B for
$m$, so $c_\ell\leq 2c_{\ell+1}$ by
Remark~\ref{rem:algoritmo}(b), and we are in~(ii).
Otherwise $m_{\ell+1}=c_{\ell+1}<c_\ell$, so $\ell$
is of type~A for $m$ and
$c_0=m_\ell=c_\ell+c_{\ell+1}$; moreover
$\ell\neq r-2$, since $m_{r-1}=c_{r-1}+1$, and we
are in~(iii). Hence every integer that can be put
in front of $\mathcal{S}$ is given by (i)--(iii).

For the converse, if $c_0\leq 2c_\ell$ and $0$ is
of type~A for some $m\in\mathcal{M}(\mathcal{S})$,
then $\ell=1$, $r\geq3$, $m_0=c_0+c_1$ and
$c_1<c_0\leq 2c_1$, so
Lemma~\ref{lem:cambio-posicion}, with $k=1$,
replaces $m_1=c_1$ with $c_0$, and $0$ is of
type~B for the new element. Hence, as
$\mathcal{M}(\mathcal{S})\neq\emptyset$, if
$c_0\leq 2c_\ell$ there is
$m\in\mathcal{M}(\mathcal{S})$ with $m_0=m_1=c_0$,
and the integers in~(i) can be put in front of
$\mathcal{S}$ through~$m$.

Finally, consider an item, (ii) or~(iii), whose
condition holds. If $0$ is of type~B for no element
of $\mathcal{M}(\mathcal{S})$, then
$c_0>2c_\ell\geq c_\ell+c_{\ell+1}$, so the item
is~(ii), and, as $0$ is of type~A for the elements
of $\mathcal{M}(\mathcal{S})$, the integer
$\ell c_0+c_\ell=c_0+c_1$ can be put in front
through them. Otherwise there is
$m\in\mathcal{M}(\mathcal{S})$ with $m_0=m_1=c_0$,
and $\ell c_0+m_{\ell+1}$ can be put in front
through $m$. If this is not the integer of our
item, we have seen that it is the integer of the
other one, whose condition also holds; thus
$\ell\leq r-3$, $c_{\ell+1}<c_\ell\leq 2c_{\ell+1}$
and $c_0=c_\ell+c_{\ell+1}$, and
Lemma~\ref{lem:cambio-posicion}, with
$k=\ell+1$, replaces $m_{\ell+1}$ with the other
element of $\{c_\ell,c_{\ell+1}\}$.
\end{proof}

\medskip

Since removing the first terms of a concentration sequence gives again a concentration sequence, namely that of a strict transform, Remark~\ref{rem:delante-conc} and Proposition~\ref{prop:imagen} decide from right to left whether a sequence $(c_0,\ldots,c_{r-1})$ of positive integers with $r\geq2$ is a concentration sequence: one checks that $c_{r-2}=c_{r-1}+1$ and then, for $i=r-3,\ldots,0$, that $c_i$ is one of the integers that Proposition~\ref{prop:imagen} allows in front of $(c_{i+1},\ldots,c_{r-1})$.
 
\begin{examples}\label{ex:primeros}
\leavevmode
\begin{enumerate}[label=(\arabic*),leftmargin=*]
\item Let us check that
$(7,5,3,2,1)$ is a concentration sequence.
First, $2=1+1$. In front of $(2,1)$, where
$\ell=1=r-1$, Proposition~\ref{prop:imagen} allows
$2$ and $3$; in front of $(3,2,1)$, where
$\ell=1=r-2$, it allows $3$ and $5$, but not
$4$; and in front of $(5,3,2,1)$ it allows $5$,
$8$ and also $7$, because $5=3+2$. On the
other hand, $(5,3,2)$ is not a concentration
sequence, since only $3$ and $4$ can be put in
front of $(3,2)$.
\item The number of integers that can be put in front of a concentration sequence $\mathcal{S}$ is at least one, since the first term of any element of $\mathcal{M}(\mathcal{S})$ can be put in front (Remark~\ref{rem:delante-conc}), and at most $\ell+2$, by Proposition~\ref{prop:imagen}. It takes each of the values $1,\ldots,5$:
\begin{itemize}
\item In front of $(5,2,2,1)$, only $7$ can be
put: here $\ell=1$ and $c_0=5>2c_1$, so (i) and
(iii) do not apply, and (ii) gives
$c_0+c_1=7$. Accordingly,
$\mathcal{M}((5,2,2,1))=\{(7,2,2,2,1)\}$
(Example~\ref{ex:arbol}).
\item In front of $(7,5,2,2,1)$, the integers
are $7$ and $9$: here $\ell=1$ and
$c_1=5>2c_2$, so (ii) does not apply, (i)
gives $7$, and (iii) gives $9$, because
$7=5+2$.
\item In front of $(5,3,2,1)$, the three items
give $5$, $8$ and $7$.
\item In front of $(3,3,3,2,1)$, where
$\ell=3=r-2$, the integers are $3$, $6$, $9$,
by~(i), and $11$, by~(ii); (iii) does not
apply, although $3=2+1$.
\item In front of $(5,5,5,3,2,1)$, where
$\ell=3$, the integers are $5$, $10$, $15$,
by~(i), $18$, by~(ii), and $17$, by~(iii), that
is, $\ell+2=5$ integers.
\end{itemize}
\end{enumerate}
\end{examples}

\section{The extended concentration sequence}
\label{sec:extendida}

In this section we show that the
conductor of the value semigroup
distinguishes the elements of
$\mathcal{M}(\mathcal{S})$
(Proposition~\ref{prop:conductor}), and
we deduce that the concentration sequence,
together with the conductor, determines the
equisingularity class of a plane algebroid
branch
(Theorem~\ref{teo:invariante-completo}).

The five branches of Example~\ref{ex:arbol}
have the same concentration sequence
$(7,4,3,2)$ but are not equisingular.
The conductors of their value semigroups
are $96$, $144$, $140$, $102$ and $170$,
respectively, which are pairwise distinct.
To compute the conductor directly from
the multiplicity sequence,
for a sequence $m=(m_0,\ldots,m_{r-1},1)$
of positive integers we set
\begin{equation}\label{eq:conductor-formula}
  \sigma(m)
  :=\sum_{k=0}^{r}m_k(m_k-1).
\end{equation}
If $m$ is the multiplicity sequence
of a plane algebroid branch $C$,
then $\sigma(m)=c(\Gamma(C))$,
the conductor of the value semigroup.
Indeed, let $\delta(C)=\dim_K\bigl(\overline{\mathcal{O}}/\mathcal{O}\bigr)$
be the $\delta$-invariant of $C$, where
$\mathcal{O}=K[[x,y]]/(f)$ and
$\overline{\mathcal{O}}$ is its integral closure;
it equals the number of elements of
$\mathbb{N}\setminus\Gamma(C)$.
Since $C$ is a plane curve,
$\mathcal{O}$ is Gorenstein,
so $c(\Gamma(C))=2\delta(C)$ (see \cite{Kunz}),
and the classical formula
$\delta(C)=\sum_{k=0}^{r}\binom{m_k}{2}$
(see \cite{Zariski-book} for $K=\mathbb{C}$
and \cite{Campillo} for $K$ of any
characteristic) gives the claim.
For instance,
$\sigma(7,7,3,3,1)=42+42+6+6=96$.
Proposition~\ref{prop:conductor}
below shows that the conductor always
distinguishes the branches with the same
concentration sequence.

The proof of
Proposition~\ref{prop:conductor}
relies on the following arithmetic lemma.

\begin{lemma}\label{lem:aritmetico}
There are no integers $w\geq 0$ and $z\geq 1$
with $w(3w-1)=z(7z-1)$.
\end{lemma}

\begin{proof}
Since $z(7z-1)>0$, we may assume $w\geq 1$.
Let $\eta=\gcd(w,z)$, $w=\eta W$, $z=\eta Z$,
so $\gcd(W,Z)=1$ and
$W(3\eta W-1)=Z(7\eta Z-1)$.
As $W$ and $Z$ are coprime,
$W$ divides $7\eta Z-1$.
Write $7\eta Z-1=uW$ with $u\geq 1$;
then $3\eta W-1=uZ$.
Multiplying the first equation by $u$
and using $uZ=3\eta W-1$,
and the second one by $u$
and using $uW=7\eta Z-1$,
we get
\[
(21\eta^2-u^2)W=u+7\eta,
\]
\[
(21\eta^2-u^2)Z=u+3\eta .
\]
Put $D:=21\eta^2-u^2$.
Since $W\geq 1$ and $u+7\eta>0$,
we get $D>0$, and $D$ divides
$(u+7\eta)-(u+3\eta)=4\eta$.

We claim that $D$ is a power of $2$.
Let $p$ be an odd prime dividing $D$.
Since $D$ divides $4\eta$, $p$ divides
$\eta$, and since $u^2=21\eta^2-D$,
$p$ divides $u$.
Moreover, if $p^a$ divides $\eta$ and $u$,
then $p^{2a}$ divides $21\eta^2-u^2=D$.
Since $D$ divides $4\eta$ and $p$ is odd,
$p^{2a}$ divides $\eta$; since $D$ divides
$u+3\eta$, $p^{2a}$ divides $u+3\eta$,
and hence $p^{2a}$ divides $u$.
Thus $p$ divides $\eta$ and $u$;
applying this with $a=1$, $p^2$ divides
$\eta$ and $u$; applying it with $a=2$,
$p^4$ divides $\eta$ and $u$; and so on.
Hence $p^{2^k}$ divides $\eta$ for every
$k\geq 0$, which is impossible, because
$\eta\geq 1$ and $p^{2^k}>\eta$
for $k$ large enough.
This proves the claim, so
$D=2^{\varepsilon}$ for some
$\varepsilon\geq 0$.

Now $D=21\eta^2-u^2=2^{\varepsilon}$, so
$u^2=21\eta^2-2^{\varepsilon}$ and, since $7$
divides $21$,
\[
u^2\equiv-2^{\varepsilon}\pmod 7.
\]
In particular $7$ does not divide $u$,
since otherwise $7$ would divide $2^{\varepsilon}$.
Hence $u^2$ is a nonzero square
modulo $7$, that is,
$u^2\equiv 1$, $2$ or $4\pmod 7$.
On the other hand,
$2^{\varepsilon}\equiv 1$, $2$ or $4\pmod 7$,
so $-2^{\varepsilon}\equiv 6$, $5$ or $3\pmod 7$,
which is a contradiction.
Hence there are no integers $w\geq 0$
and $z\geq 1$ with $w(3w-1)=z(7z-1)$.
\end{proof}

\begin{remark}\label{rem:aritmetica}
Lemma~\ref{lem:aritmetico} is essential
in the proof of
Proposition~\ref{prop:conductor}:
for the pairs in the fifth row of
Table~\ref{tab:pares}, which correspond
to case~(III) of
Lemma~\ref{lem:tramo-inicial}
(see Remark~\ref{rem:tabla}),
the difference of the conductors of two
branches with the same concentration
sequence is $z(7z-1)-w(3w-1)$, which can
be positive or negative.
For instance, the multiplicity sequences
$(10,7,3,3,1)$ and $(11,4,4,3,1)$
of Example~\ref{ex:arbol} correspond to
$z=1$ and $w=2$, and their conductors
$144$ and $140$ differ only by $4$.
Note also that the similar equation
$w(3w+1)=z(7z-1)$ does have solutions,
for instance $w=15$ and $z=10$.
\end{remark}

To prove that the conductor distinguishes
the sequences with the same concentration
sequence, we look at the last position where
two such sequences differ.
We write $\psi(x)=x(x-1)$, so that
$\sigma(m)=\sum_{k=0}^{r}\psi(m_k)$, and we
use that $\psi(x)-\psi(y)=(x-y)(x+y-1)$.

\begin{lemma}\label{lem:ultima}
Let $\mathcal{S}=(c_0,\ldots,c_{r-1})$ be a
sequence of positive integers, let
$m=(m_0,\ldots,m_{r-1},1)$ and
$\tilde m=(\tilde m_0,\ldots,\tilde m_{r-1},1)$
be two elements of
$\mathcal{M}(\mathcal{S})$, assume
that $m\neq\tilde m$, and let $i$ be the
largest index with $m_i\neq\tilde m_i$.
Then $i\leq r-2$, $c_i=m_{i+1}$ and, after
exchanging $m$ and $\tilde m$ if necessary,
$m_i=m_{i+1}+m_{i+2}$ and
$\tilde m_i=m_{i+1}$. Writing
\[
z=m_{i+2},\qquad w=m_{i+1}-m_{i+2},
\qquad \Lambda=2w+3z,
\]
we have $z\geq 1$, $w\geq 0$,
$m_i=w+2z$, $\tilde m_i=w+z$,
$m_i+\tilde m_i=\Lambda$, and
\begin{equation}\label{eq:D}
D:=\sigma(m)-\sigma(\tilde m)
=z(\Lambda-1)
+\sum_{k<i}\bigl(\psi(m_k)-\psi(\tilde m_k)\bigr).
\end{equation}
\end{lemma}

\begin{proof}
By Theorem~\ref{teo:descripcion}, applied
to $m$ and $\tilde m$, which belong to
$\mathcal{M}(\mathcal{S})$, we have
$m_{r-1}=\tilde m_{r-1}=c_{r-1}+1$; since
also $m_r=\tilde m_r=1$, we get $i\leq r-2$.
Since $i$ is the largest index where $m$ and
$\tilde m$ differ, they agree at the
positions $i+1$ and $i+2$, that is,
$m_{i+1}=\tilde m_{i+1}$ and
$m_{i+2}=\tilde m_{i+2}$.

We now show that the index $i$ is of type~B
for both sequences. Recall that $i$ is of
type~A for $m$ (respectively $\tilde m$) if
and only if $c_i>m_{i+1}$ (respectively
$c_i>\tilde m_{i+1}$). Since $c_i$ is the
same for both sequences and
$m_{i+1}=\tilde m_{i+1}$, these two
inequalities are the same; hence $i$ is
either of type~A for both sequences or of
type~B for both.
If it were of type~A, case~(A) of
Theorem~\ref{teo:descripcion} would give
$m_i=c_i+m_{i+1}=c_i+\tilde m_{i+1}
=\tilde m_i$, which is false. Hence $i$ is of
type~B for both sequences, that is,
$c_i=m_{i+1}$, and case~(B) of
Theorem~\ref{teo:descripcion} gives
$m_i,\tilde m_i\in
\{m_{i+1},m_{i+1}+m_{i+2}\}$ (recall that
$m$ and $\tilde m$ agree at $i+1$ and $i+2$).
Since $m_i\neq\tilde m_i$ and both belong to
this set, one of them is $m_{i+1}$ and the
other is $m_{i+1}+m_{i+2}$; exchanging $m$
and $\tilde m$ if necessary, we may assume
that $m_i=m_{i+1}+m_{i+2}$ and
$\tilde m_i=m_{i+1}$.

Now $z\geq 1$, because it is a
multiplicity, and $w\geq 0$, because $m$ is
non-increasing. As $m_{i+1}=w+z$, we get
$m_i=m_{i+1}+m_{i+2}=w+2z$,
$\tilde m_i=m_{i+1}=w+z$ and
$m_i+\tilde m_i=2w+3z=\Lambda$.

Finally, by the definition of $\sigma$,
$D=\sum_{k=0}^{r}\bigl(\psi(m_k)
-\psi(\tilde m_k)\bigr)$. The terms with
$k>i$ vanish, because $m_k=\tilde m_k$ for
$k>i$, and, by the identity
$\psi(x)-\psi(y)=(x-y)(x+y-1)$ applied with
$x=m_i=w+2z$ and $y=\tilde m_i=w+z$, the
term with $k=i$ is
\[
\psi(m_i)-\psi(\tilde m_i)
=z\,(2w+3z-1)=z(\Lambda-1).
\]
This gives~\eqref{eq:D}.
\end{proof}

\begin{lemma}\label{lem:casos}
Let $\mathcal{S}$, $m$, $\tilde m$, $i$, $z$,
$w$, $\Lambda$ and $D$ be as in
Lemma~\ref{lem:ultima}, and assume that
$i\leq 2$. Then $(m_0,\ldots,m_i)$ and
$(\tilde m_0,\ldots,\tilde m_i)$ are given by
one of the nine rows of
Table~\ref{tab:pares}, $D$ is the value
given there, and $D\neq 0$.
\end{lemma}

\begin{table}[ht]
\centering
\renewcommand{\arraystretch}{1.15}
\begin{tabular}{cllr}
\toprule
$i$ & $(m_0,\ldots,m_i)$ & $(\tilde m_0,\ldots,\tilde m_i)$ & $D$\\
\midrule
$0$ & $(w+2z)$ & $(w+z)$ & $z(2w+3z-1)$\\
\midrule
$1$ & $(w+2z,\,w+2z)$ & $(\Lambda,\,w+z)$ & $-3w^2-6wz-2z^2+w$\\
$1$ & $(\Lambda,\,w+2z)$ & $(\Lambda,\,w+z)$ & $z(2w+3z-1)$\\
\midrule
$2$ & $(3w+5z,\,w+2z,\,w+2z)$ & $(\Lambda,\,\Lambda,\,w+z)$ & $2w^2+12wz+14z^2-2z$\\
$2$ & $(3w+5z,\,w+2z,\,w+2z)$ & $(3w+4z,\,\Lambda,\,w+z)$ & $z(7z-1)-w(3w-1)$\\
$2$ & $(3w+5z,\,\Lambda,\,w+2z)$ & $(3w+4z,\,\Lambda,\,w+z)$ & $2z(4w+6z-1)$\\
$2$ & $(3w+5z,\,\Lambda,\,w+2z)$ & $(\Lambda,\,\Lambda,\,w+z)$ & $5w^2+20wz+19z^2-w-3z$\\
$2$ & $(\Lambda,\,\Lambda,\,w+2z)$ & $(3w+4z,\,\Lambda,\,w+z)$ & $-5w^2-10wz-4z^2+w$\\
$2$ & $(\Lambda,\,\Lambda,\,w+2z)$ & $(\Lambda,\,\Lambda,\,w+z)$ & $z(2w+3z-1)$\\
\bottomrule
\end{tabular}
\caption{The pairs in Lemma~\ref{lem:casos}.
Here $\Lambda=2w+3z$, so that
$3w+5z=\Lambda+w+2z$ and
$3w+4z=\Lambda+w+z$.}\label{tab:pares}
\end{table}

\begin{proof}
From Lemma~\ref{lem:ultima} we have
$c_i=m_{i+1}$, $m_i=w+2z$ and
$\tilde m_i=w+z$; moreover,
$m_{i+1}=\tilde m_{i+1}=w+z$ by the
definitions of $w$ and $z$. We first determine the
terms at the positions $i-1$ and $i-2$.

Assume first that $i\in\{1,2\}$. We
determine $c_{i-1}$, $m_{i-1}$ and
$\tilde m_{i-1}$.
\begin{itemize}
\item In $m$: we have
$c_{i-1}=\max\{m_i,m_{i-1}-m_i\}\geq m_i$ by
Proposition~\ref{proposicion concentracion curva}.
If $c_{i-1}>m_i$, then
Lemma~\ref{lem:tipoB}(b), applied with
$i-1$ in place of $i$, yields
$m_{i+1}=m_i$, which is false, since
$m_{i+1}=w+z<w+2z=m_i$. Hence
$c_{i-1}=m_i=w+2z$, and, by
Lemma~\ref{lem:tipoB}(a), applied with
$i-1$ in place of $i$,
$m_{i-1}\in\{m_i,m_i+m_{i+1}\}$, that is,
$m_{i-1}\in\{w+2z,\Lambda\}$, since
$m_i+m_{i+1}=2w+3z=\Lambda$.
\item In $\tilde m$: since
$c_{i-1}=w+2z>w+z=\tilde m_i$, case~(A) of
Theorem~\ref{teo:descripcion} holds at
$i-1$; thus
$\tilde m_{i-1}=c_{i-1}+\tilde m_i
=(w+2z)+(w+z)=\Lambda$.
\end{itemize}
Assume moreover that $i=2$. We determine
$c_{i-2}$, $m_{i-2}$ and $\tilde m_{i-2}$,
that is, $c_0$, $m_0$ and $\tilde m_0$.
\begin{itemize}
\item In $\tilde m$: by
Proposition~\ref{proposicion concentracion curva},
$c_{i-2}\geq\tilde m_{i-1}$. If
$c_{i-2}>\tilde m_{i-1}$, then
Lemma~\ref{lem:tipoB}(b), applied with
$i-2$ in place of $i$, yields
$\tilde m_i=\tilde m_{i-1}$, which is false,
since $\tilde m_i=w+z<\Lambda=\tilde m_{i-1}$.
Hence $c_{i-2}=\tilde m_{i-1}=\Lambda$, and,
by Lemma~\ref{lem:tipoB}(a),
$\tilde m_{i-2}\in\{\Lambda,\Lambda+(w+z)\}
=\{\Lambda,3w+4z\}$.
\item In $m$: recall that
$m_{i-1}\in\{w+2z,\Lambda\}$ and
$c_{i-2}=\Lambda$. If $m_{i-1}=w+2z$, then
$c_{i-2}=\Lambda>m_{i-1}$, so case~(A) of
Theorem~\ref{teo:descripcion} holds at
$i-2$; thus
$m_{i-2}=c_{i-2}+m_{i-1}=\Lambda+(w+2z)
=3w+5z$. If $m_{i-1}=\Lambda$, then
$c_{i-2}=\Lambda=m_{i-1}$, so case~(B)
holds at $i-2$, and
$m_{i-2}\in\{\Lambda,\Lambda+m_i\}
=\{\Lambda,3w+5z\}$.
\end{itemize}
Since $i\leq 2$, this determines all the
possibilities for $(m_0,\ldots,m_i)$ and
$(\tilde m_0,\ldots,\tilde m_i)$: for $i=0$
we get the first row of
Table~\ref{tab:pares}; for $i=1$, the two
choices of $m_0$ give the second and third
rows; and for $i=2$, the three possibilities
for $(m_0,m_1)$ and the two for
$(\tilde m_0,\tilde m_1)$ give the remaining
six rows.

We now compute $D$. By~\eqref{eq:D},
$D$ is $z(\Lambda-1)=z(2w+3z-1)$ plus the
differences $\psi(m_k)-\psi(\tilde m_k)$ at
the positions $k<i$ where $m$ and
$\tilde m$ differ. By the identity
$\psi(x)-\psi(y)=(x-y)(x+y-1)$, these
differences are given by
\[
\begin{aligned}
\psi(\Lambda)-\psi(w+2z)&=(w+z)(3w+5z-1),\\
\psi(3w+5z)-\psi(\Lambda)&=(w+2z)(5w+8z-1),\\
\psi(3w+5z)-\psi(3w+4z)&=z(6w+9z-1),\\
\psi(3w+4z)-\psi(\Lambda)&=(w+z)(5w+7z-1).
\end{aligned}
\]
For instance, in the fifth row the
sequences differ at the positions $0$ and
$1$, and
\[
\begin{aligned}
D&=z(2w+3z-1)+z(6w+9z-1)
-(w+z)(3w+5z-1)\\
&=z(7z-1)-w(3w-1).
\end{aligned}
\]
The other rows are computed in the same way,
and give the values of
Table~\ref{tab:pares}.

Finally, we show that $D\neq 0$, using that
$z\geq 1$ and $w\geq 0$.
\begin{itemize}
\item In the first, third and ninth rows,
$D=z(2w+3z-1)>0$.
\item In the fourth, sixth and seventh
rows, $D>0$, since
$2w^2+12wz+14z^2-2z\geq 14z^2-2z>0$,
$2z(4w+6z-1)>0$ and
$5w^2+20wz+19z^2-w-3z
=w(5w-1)+20wz+z(19z-3)>0$.
\item In the second and eighth rows, $D<0$,
since $D=-w(3w+6z-1)-2z^2$ and
$D=-w(5w+10z-1)-4z^2$, respectively.
\item In the fifth row,
$D=z(7z-1)-w(3w-1)\neq 0$ by
Lemma~\ref{lem:aritmetico}.
\end{itemize}
\end{proof}

\begin{remark}\label{rem:tabla}
The rows of Table~\ref{tab:pares}
correspond to the situations described in
Lemmas~\ref{lem:dif-aislada}
and~\ref{lem:tramo-inicial}.
The first row is case~(I) of
Lemma~\ref{lem:tramo-inicial}, with no
further difference; the second row is
case~(II); the fourth and fifth rows are
case~(III), with $m$ and $\tilde m$
exchanged, and with $\tilde m_0=c_0$ and
$\tilde m_0=c_0+\tilde m_2$, respectively;
the sixth, seventh and eighth rows are
case~(I) together with a difference at the
position~$2$, which is the case excluded in
Lemma~\ref{lem:divisibilidad} for $k=2$;
and in the third and ninth rows
$m_0=\tilde m_0$, and the only difference is
at the position $1$ or $2$, as in
Lemma~\ref{lem:dif-aislada}(a).
\end{remark}

\begin{lemma}\label{lem:induccion}
Let $\mathcal{S}$, $m$, $\tilde m$, $i$, $z$,
$w$, $\Lambda$ and $D$ be as in
Lemma~\ref{lem:ultima}, and assume that
$i\geq 3$. Then $m_{i-1}=\tilde m_{i-1}=\Lambda$, all the $m_k$ and $\tilde m_k$ with
$k\leq i-1$ are multiples of $\Lambda$, and
$D\equiv -z\pmod\Lambda$.
In particular $D\neq 0$.
\end{lemma}

\begin{proof}
Let $t$ be the least index with
$m_t=\tilde m_t$. If $m_0=\tilde m_0$, then
$t=0$; otherwise $t\leq 3$ by
Lemma~\ref{lem:tramo-inicial}.
Since $t\leq 3\leq i$ and $i\neq t$ (because
$m_i\neq\tilde m_i$), we get $i>t$, so $m$ and $\tilde m$ agree at $i-1$ by Lemma~\ref{lem:dif-aislada}(b), and Lemma~\ref{lem:dif-aislada}(a), applied with $k=i$, gives $m_{i-1}=\tilde m_{i-1}=m_i+\tilde m_i=\Lambda$.
By Lemma~\ref{lem:divisibilidad}(b),
applied with $k=i\geq 3$, $\Lambda$ divides
$m_j$ and $\tilde m_j$ for every
$j\leq i-1$, and hence also
$\psi(m_j)=m_j(m_j-1)$ and
$\psi(\tilde m_j)$. So every summand of the
sum in~\eqref{eq:D} is a multiple of
$\Lambda$, and
$D\equiv z(\Lambda-1)\equiv -z\pmod\Lambda$.
Since $0<z<2w+3z=\Lambda$, the integer $-z$
is not a multiple of $\Lambda$, and hence
$D\neq 0$.
\end{proof}

\begin{proposition}\label{prop:conductor}
Let $\mathcal{S}=(c_0,\ldots,c_{r-1})$ be a
sequence of positive integers, and let $m$
and $\tilde m$ be two elements of
$\mathcal{M}(\mathcal{S})$.
If $\sigma(m)=\sigma(\tilde m)$,
then $m=\tilde m$.
\end{proposition}

\begin{proof}
Suppose that $m\neq\tilde m$, and let $i$,
$z$, $w$, $\Lambda$ and
$D=\sigma(m)-\sigma(\tilde m)$ be as in
Lemma~\ref{lem:ultima}.
By Lemma~\ref{lem:casos} if $i\leq 2$, and by
Lemma~\ref{lem:induccion} if $i\geq 3$, we
get $D\neq 0$, that is,
$\sigma(m)\neq\sigma(\tilde m)$.
\end{proof}

\begin{remark}\label{rem:hipotesis-rama}
The hypothesis that $m$ and $\tilde m$ are
multiplicity sequences of branches cannot be
dropped in Proposition~\ref{prop:conductor}.
It is used in the proof through Theorem~\ref{teo:descripcion} and
Lemmas~\ref{lem:tipoB},
\ref{lem:dif-aislada},
\ref{lem:tramo-inicial}
and~\ref{lem:divisibilidad}, whose
proofs rely on conditions~(1) and~(2) of
Proposition~\ref{prop:mult-to-semigrup},
which characterize multiplicity sequences
of branches.
For instance, the non-increasing sequences
$(16,7,7,1)$ and $(15,9,7,1)$ give the same
sequence
$\bigl(\max\{m_{i+1},m_i-m_{i+1}\}\bigr)_i=(9,7,6)$
and the same value $\sigma=324$, but neither
satisfies condition~(1): $7$ does not divide
$16$ and $16-2\cdot 7=2\neq 1$, and $9$ does
not divide $15$ and $15-9=6\neq 7$.
\end{remark}

By Proposition~\ref{prop:conductor},
the concentration sequence of a plane
algebroid branch, together with the
conductor of its value semigroup, determines its multiplicity sequence. This motivates the following definition.

\begin{definition}\label{def:extended-concentration}
Let $C$ be a plane algebroid branch,
let $r$ be the number of blowups needed
to reach the first smooth strict transform
of $C$, and let
$\mathcal{S}=(c_0,\ldots,c_{r-1})$
be its concentration sequence.
The \emph{extended concentration sequence}
of $C$ is
\[
\mathcal{S}^e(C)=(c(\Gamma(C));\,c_0,\ldots,c_{r-1}),
\]
where $c(\Gamma(C))$ is the conductor
of $\Gamma(C)$.
\end{definition}

\begin{theorem}\label{teo:invariante-completo}
Two plane algebroid branches are
equisingular if and only if they have
the same extended concentration sequence.
\end{theorem}

\begin{proof}
Equisingular branches have the same value semigroup, hence the same multiplicity sequence by Proposition~\ref{prop:Angermuller}, and therefore the same concentration sequence and the same conductor. Conversely, let $C$ and $\tilde C$ be branches with the same extended concentration sequence. In particular they have the same concentration sequence $\mathcal{S}=(c_0,\ldots,c_{r-1})$.
Since $\mathcal{S}$ has $r$ terms, both branches need $r$ blowups to reach a smooth strict transform, so their multiplicity sequences $m=(m_0,\ldots,m_{r-1},1)$ and $\tilde m=(\tilde m_0,\ldots,\tilde m_{r-1},1)$ have the same length.
Moreover, $\sigma(m)=c(\Gamma(C))=c(\Gamma(\tilde C))=\sigma(\tilde m)$,
so we can apply Proposition~\ref{prop:conductor}, which gives $m=\tilde m$. By Proposition~\ref{prop:mult-to-semigrup}, $C$ and $\tilde C$ have the same value semigroup, so they are equisingular.
\end{proof}

\begin{remark}\label{rem:milnor}
Neither of the two data in Definition~\ref{def:extended-concentration}
is enough by itself to determine the equisingularity class of a plane algebroid branch. The branches of Example~\ref{ex:arbol} have the same concentration sequence $(7,4,3,2)$ and are not equisingular. On the other hand, the branches with value semigroups $\langle 4,6,25\rangle$ and $\langle 4,10,21\rangle$ have the same multiplicity, the same genus and the same conductor $28$, but they are not equisingular, since their value semigroups are different. Their concentration sequences, $(2,2,2,2,2,2,2,2,1)$ and $(4,2,2,1)$, respectively, are different, so in this example it is the concentration sequence that distinguishes them.

If $K$ has characteristic zero, then the conductor of $\Gamma(C)$ equals the Milnor number of $C$, since both are equal to $2\delta(C)$ (see~\cite[Theorem~10.5]{Milnor} and~\cite{Kunz}). Hence, in characteristic zero, Theorem~\ref{teo:invariante-completo} says that the concentration sequence and the Milnor number together form a complete
invariant of the equisingularity class,
while neither of them is so by itself. In positive characteristic this equality may fail, since the Milnor number can be larger than the conductor (see~\cite{Melle-Wall}).
\end{remark}

\medskip
\noindent\textbf{Funding.} The authors were
supported by the Spanish grant
PID2023-149508NB-I00, funded by
MICIU/AEI/10.13039/501100011033 and by ERDF/EU.

\medskip
\noindent\textbf{Acknowledgements.} The authors would like to thank Marco D'Anna for a question that led to the development of some of the results in this paper. We also thank Marcelo E. Hernandes for his comments on earlier versions of the manuscript.

\noindent {\sc Evelia R. Garc\'ia Barroso}\\
Departamento de Matem\'aticas, Estad\'{\i}stica e Investigaci\'on Operativa\\
	Instituto Universitario de Matem\'aticas y Aplicaciones (IMAULL)\\ Universidad de La Laguna. Apartado de Correos 456.
	38200 La Laguna, Tenerife, Espa\~na \\ORCID ID: 0000-0001-7575-2619
	
	 \noindent {ergarcia@ull.es}
	\vspace{0.3cm}

\noindent {\sc William G. Hern\'andez-Yanes}\\
Departamento de Matem\'aticas, Estad\'{\i}stica e Investigaci\'on Operativa\\
	Instituto Universitario de Matem\'aticas y Aplicaciones (IMAULL)\\ Universidad de La Laguna. Apartado de Correos 456.
	38200 La Laguna, Tenerife, Espa\~na \\ORCID ID: 0009-0002-5140-8295
	
	 \noindent {whernany@ull.es}
	\vspace{0.3cm}

\end{document}